\documentclass[reqno, 11pt, a4paper, final]{amsart}
\usepackage{amssymb,amsmath,amsthm}
\usepackage[foot]{amsaddr}

\usepackage[english]{babel}
\usepackage[T1]{fontenc}
\usepackage[utf8]{inputenc}

\usepackage{hyperref,url,doi}
\usepackage{graphicx}
\usepackage{tikz-cd}
\usepackage{subcaption}
\usepackage{emptypage}
\usepackage{faktor}
\usepackage{enumitem}
\usepackage{mathrsfs}
\usepackage{ragged2e}
\usepackage[
top=2.5cm,
bottom=2.5cm,
inner=3cm, 
outer=3cm,
footskip=40pt 
]{geometry}

\usepackage{physics}

\numberwithin{equation}{subsection}

\theoremstyle{definition}
\newtheorem{defi}{Definition}[section]

\newtheorem{rmk}[defi]{Remark}

\theoremstyle{plain}
\newtheorem{lm}[defi]{Lemma}
\newtheorem{thm}[defi]{Theorem}

\newtheorem{cor}[defi]{Corollary}

\theoremstyle{remark}
\newtheorem{claim}[defi]{Claim}

\newcommand{\R}{\mathbb R}

\newcommand{\N}{\mathbb N}
\newcommand{\Z}{\mathbb Z}
\newcommand{\C}{\mathbb C}

\newcommand{\D}{\mathbb D}
\newcommand{\T}{\mathbb T}
\newcommand{\bx}{\mathbf{x}}
\newcommand{\by}{\mathbf{y}}

\renewcommand{\c}{\mathcal{C}}
\renewcommand{\r}{\mathcal{R}}
\renewcommand{\a}{\mathcal{A}}
\renewcommand{\o}{\mathcal{O}}

\newcommand{\coD}{\mathscr D}
\newcommand{\coL}{\mathscr L}

\DeclareMathOperator{\dist}{dist}
\DeclareMathOperator{\pr}{pr}
\DeclareMathOperator{\Int}{Int}
\newcommand{\mmin}{\text{min}}
\newcommand{\mmax}{\text{max}}

\usepackage{lmodern}
\newlength\longest

\begin{document}
	
	\title{Holomorphic Legendrian curves in convex domains}
	\author[Andrej Svetina]{Andrej Svetina}
	
	
	\address[Andrej Svetina]{Faculty of Mathematics and Physics, University of Ljubljana, Jadranska ulica 21, 1000 Ljubljana, Slovenia}
	\email{andrej.svetina@fmf.uni-lj.si}

	
	\date{\today}
	

	\begin{abstract}
	We prove several results on approximation and interpolation of holomorphic Legendrian curves in convex domains in $\C^{2n+1}$, $n \geq 2$, with the standard contact structure. Namely, we show that such a curve, defined on a compact bordered Riemann surface $M$, whose image lies in the interior of a convex domain $\coD \subset \C^{2n+1}$, may be approximated uniformly on compacts in the interior $\Int M$ by holomorphic Legendrian curves $\Int M \to \coD$ such that the approximants are proper, complete, agree with the starting curve on a given finite set in $\Int M$ to a given finite order, and hit a specified diverging discrete set in the convex domain. We first show approximation of this kind on bounded strongly convex domains and then generalise it to arbitrary convex domains. As a consequence we show that any bordered Riemann surface properly embeds into a convex domain as a complete holomorphic Legendrian curve under a suitable geometric condition on the boundary of the codomain.
\end{abstract}

\maketitle
	
	\section{Introduction}
	The standard complex contact structure on $\C^{2n+1}$ is the holomorphic hyperplane subbundle in the holomorphic tangent bundle, given as the kernel of the holomorphic differential form
	\[
	\alpha = \dd z + x_1 \dd y_1 + \cdots x_n \dd y_n,
	\]
	where $(x_1,y_1, \ldots, x_n,y_n,z) \in \C^{2n+1}$ are standard complex coordinates on $\C^{2n+1}$. Note that the subbundle defined by $\xi = \ker \alpha$ is \emph{completely nonintegrable}, meaning that the wedge product $\alpha \wedge (\dd \alpha)^n$ defines the standard complex Euclidean volume form at every point $p \in \C^{2n+1}$, thus the condition $\alpha \wedge \dd \alpha = 0$ of the Frobenius theorem is violated, in a sense, as much as possible. It follows in particular that the maximal complex dimension of an integral complex manifold of $\xi$ is $n$. Of particular interest are complex integral curves of $\xi$, called \emph{holomorphic Legendrian curves}, since they are the integral manifolds that are most easily manipulated. Namely, a holomorphic Legendrian curve is a holomorphic map $f\colon \r \to \C^{2n+1}$, defined on a Riemann surface $\r$, satisfying $f^*\alpha = 0$ at every point $p \in \r$. The Riemann surface in question can either be open or \emph{compact bordered}, where the latter denotes a compact surface with boundary that is conformally diffeomorphic to the closure of a domain in an open Riemann surface and whose boundary consists of finitely many Jordan curves, see \cite[Definition 1.10.8.]{mincplx}. The interior of such a surface is called a \emph{bordered Riemann surface}. It was proven by Stout in \cite{Stout1965} that a compact bordered Riemann surface is conformally diffeomorphic to a complement of a disjoint union of finitely many disks in a compact Riemann surface.
	
	It was proven by Alarcón, Forstnerič and López in \cite{Alarcon2017} that Legendrian curves in $\C^{2n+1}$ admit analogues of the classical Runge and Mergelyan approximation theorems. More precisely, if a map $f\colon K \to \C^{2n+1}$, defined on an $\o(\r)$-convex compact set $K \subset \r$ in an open Riemann surface $\r$ satisfies $f^*\alpha = 0$ on a neighbourhood of $K$, it may be approximated uniformly on $K$ by entire Legendrian curves $F \colon \r \to \C^{2n+1}$. Furthermore, they showed the approximants could be taken to satisfy certain global properties, such as embeddedness, properness, and completeness, where the latter means that the Riemannian metric $f^*g$, obtained as the pullback of the standard Euclidean metric $g$ on $\C^{2n+1}$, is complete, see Section \ref{cpltMetric}.
	
	It is notable that the authors' methods in \cite{Alarcon2017} mirror the ones they, together with D.\ Drnovšek, used in \cite{Alarcon2015} to obtain similar results for minimal surfaces. In the latter paper they have shown in addition that minimal surfaces defined on compact bordered Riemann surfaces may be approximated with proper and complete ones not only when the ambient space is, in this case real, Euclidean space, but also when the image of the curve lies in a bounded convex domain in the said Euclidean space. In the present paper we show that a similar result holds for holomorphic Legendrian curves, where by a holomorphic Legendrian curve we mean a continuously differentiable map $f\colon M \to \C^{2n+1}$, defined on a compact bordered Riemann surface $M$, that is holomorphic on $\Int M$ and satisfies $f^*\alpha = 0$ along $M$.
	
	\begin{thm}
		\label{properCvxBasic}
		Suppose $\coD$ is a smoothly bounded strongly convex domain in $\C^{2n+1}$, $n \geq 2$, and suppose $f\colon M \to \coD$ is a Legendrian curve, defined on a compact bordered Riemann surface $M$. Then, $f$ may be approximated uniformly on compacts in $\Int M$ by continuous maps $F\colon M \to \overline{\coD}$ such that $F|_{\Int M}\colon \Int M \to \coD$ is a proper and complete Legendrian embedding.
	\end{thm}
	
	The restriction $n \geq  2$ is explained in Section \ref{propEmbHyperplane} and the proof of Lemma \ref{cvxInductiveLemma}. When $2n+1$ is the dimension of a complex Euclidean space, we assume $n \geq 2$ throughout the paper, unless stated otherwise. Theorem \ref{properCvxBasic} is proven in Section \ref{pfsMain} as a consequence of the more general Theorem \ref{properCvx} which also enables one to interpolate the starting curve to a given finite order in a given finite set in $\Int M$. Additionally one is able to make the approximating curve hit a specified closed and discrete set in the codomain, however at the cost of losing continuity up to the boundary. On the other hand, dropping this assumption also enables approximation in more general convex domains such as weakly convex, unbounded domains or domains whose boundary is not necessarily smooth. Similarly as in \cite{Alarcon2015}, Theorem \ref{weakCvx} is proven using an inductive procedure over an exhaustion of the domain by strongly convex ones. The precise result is stated in Theorem \ref{convExhaustionThm}, below we only state some interesting consequences.

	\begin{thm}
		\label{weakCvx}
		Suppose $\coD$ is a convex domain in $\C^{2n+1}$, $n \geq 2$. and $f\colon M \to \coD$ is a holomorphic Legendrian curve, defined on a compact bordered Riemann surface $M$. For a given closed discrete set $\Lambda \subset \coD$ the curve $f$ may be approximated uniformly on compacts in $\Int M$ by proper and complete holomorphic Legendrian embeddings $F\colon \Int M \to \coD$ such that $\Lambda \subset F(\Int M)$.
	\end{thm}
	
	This result follows from the more precise Theorem \ref{convExhaustionThm}. Another of its consequences is the following existence and hitting theorem for bordered Riemann surfaces.

	\begin{cor}
		\label{properExst}
		Suppose $\coD$ is a convex domain in $\C^{2n+1}$, $n \geq 2$, $\Lambda \subset \coD$ is a closed discrete subset and $M$ is a bordered Riemann surface. There exists a proper and complete holomorphic Legendrian embedding $F\colon M \to \coD$ hitting every point $p \in \Lambda$.
	\end{cor}
	
	In particular, the unit disk $\D$ may be properly embedded into any bounded strongly convex domain $\overline{\coD} \subset \C^{2n+1}$ as a complete holomorphic Legendrian curve such that the boundary circle $b\D$ accumulates everywhere on the boundary $b\coD$. If we drop the assumption of embeddedness, we are able to make a Legendrian curve hit an arbitrary countable set in a given convex domain. A closely related result was obtained by Alarcón and Forstnerič in \cite[Theorem 10.4]{Alarcon2023}, namely, that such curves exist in arbitrary complex contact manifolds. By restricting ourselves to only convex domains in Euclidean spaces, we are able to ensure in addition that the Legendrian curve in question is \emph{almost proper}, that is, preimages of compact sets have relatively compact connected components.
	
	\begin{thm}
		\label{denseImmersion}
		Any bordered Riemann surface $M$ can be mapped to any given convex domain $\coD$ in $\C^{2n+1}$, $n \geq 2$, as an injectively immersed almost proper holomorphic Legendrian curve with dense image.
	\end{thm}
	
	This is proven in Section \ref{almostProperProof} as a special consequence of Theorem \ref{almostProperApprox} on approximation of Legendrian curves by almost proper ones.
	
	In \cite{Alarcon2015}, the authors show in addition that uniform approximation of already proper minimal surfaces with proper and complete ones is possible. The key step in the proof consists of using a translation and a homothety to perturb the proper curve into one whose image lies entirely in the interior of the convex domain. This is possible by the fact that such transformations map minimal surfaces into minimal surfaces. In Section \ref{properWcomplete} we outline one possible analogue of this theorem by adapting translations and homotheties to the contact setting. Using such a method however requires an additional geometric assumption on the boundary of the convex domain in question.

	\section{Preliminaries}
	
	In this section we collect some technical results that will be needed in the proof.
	
	\subsection{Properly embedded Legendrian lines in affine hyperplanes in $\C^{2n+1}$}
	\label{propEmbHyperplane}
	
	A crucial step in constructing proper Legendrian curves in a convex domain $\coD \subset \C^{2n+1}$ involves finding certain complex Legendrian lines $\Psi \colon \C \to \C^{2n+1}$ lying in a given complex hyperplane $\Pi$ in $\C^{2n+1}$. One then intersects the image $\Psi(\C)$ with the convex domain $\coD$, and, provided $\Psi(C)$ intersects the boundary $b\coD$ transversely, obtains a Legendrian curve $\Psi\colon \C \supset D \to \Pi$, where $D$ is biholomorphic to the unit disk $\D \subset \C$. In fact, one wishes to find a family $\Psi_t \colon D_t \to \Pi_t$ of such Legendrian disks for $t \in [0,1]$, varying continuously in $t$, where $\Pi_t$ is a continuous family of hyperplanes in $\C^{2n+1}$ and $D_t$ is the corresponding family of proper simply connected domains $D_t \subset \C$ such that $\Psi_t$ maps $D_t$ biholomorphically onto $\Psi_t(\C) \cap \coD$. If $n \geq 2$, then there exist linear Legendrian disks in every complex hyperplane $\Pi \subset \C^{2n+1}$, so that the above intersection is always transverse. If $n=1$ however, the intersection of $\Pi$ with the contact structure in a generic point $p \in \Pi$ forms a complex one dimensional subspace in $T_p\C^3$, thus there is in general a unique Legendrian curve passing through $p$ and tangent to both $\Pi$ and $\xi = \ker \alpha$. Since we cannot guarantee that the intersection of this curve with the boundary of an arbitrary convex domain is transverse, this construction does not furnish the desired family of disks.
	
	Label the coordinates on $\C^{2n+1}$ by $(x_1, \ldots, x_n, y_1, \ldots, y_n, z)$ and let $(a_1,b_1,\ldots,a_n, b_n,) \in \C^{2n}$ be a vector, satisfying
	$$
	a_1 \cdots a_n \neq 0.
	$$
	Denote by $\Pi$ the complex hyperplane
	$$
	\Pi = \left\{(x_1, \ldots, x_n, y_1, \ldots, y_n, z) \in \C^{2n+1} \, : \, z + \sum_{j=1}^n (a_jx_j + b_jy_j) = 0\right\}
	$$
	with complex normal vector $\vec{n} = (a_1,b_1,\ldots,a_n,b_n,1) \in \C^{2n+1}$ (with respect to the standard Hermitian product on $\C^{2n+1}$). Assume first $n = 2$ and let $p_0 = (x_{1,0},x_{2,0},y_{1,0},y_{2,0}z_0) \in \C^5$ be any point. Define holomorphic functions
	\[
	X_1(\zeta) = a_2\zeta + x_{1,0}, \quad 
	X_2(\zeta) = -a_1\zeta + x_{2,0}, \quad
	Y_1 \equiv y_{1,0}, \quad
	Y_2 \equiv y_{2,0}, \quad
	Z \equiv z_0.
	\]
	Then the mapping $\Psi_{\vec{n}} \colon \C \to \C^5$, $\Psi_{\vec{n}} = (X_1,X_2,Y_1,Y_2,Z)$, defines a linear holomorphic Legendrian embedding whose image lies in the hyperplane $p_0 + \Pi$. For $n \geq 3$ such linear embeddings are obtained by an analogous construction.
	
	Suppose $n=1$ and $\Psi \colon \C \to \C^3$ is a holomorphic Legendrian curve lying in the hyperplane $\Pi \subset \C^3$ defined by the equation $z+ax+by=0$. Taking the derivative together with the contact condition implies
	\[
	0=\dd z + a\dd x+ b\dd y = -x \dd y + a \dd x + b \dd y, \quad \dv{x}{y} = \frac{x-b}{a},
	\]
	provided $a \neq 0$. The above differential equation for $x=x(y)$ has a unique solution
	\[
	x(y) = (x_0 -b) \exp((y-y_0)/a) +b.
	\]
	Choosing the parametrisation $y = \zeta + y_0$ we obtain
	\begin{align*}
		X(\zeta) &= (x_0 -b) \exp(\zeta/a) +b\\
		Y(\zeta) &= \zeta + y_0 \\
		Z(\zeta) &= z_0 + ax_0 + by_0 - aX(\zeta) - bY(\zeta),
	\end{align*}
	and label the corresponding Legendrian curve $\C \to \C^3$ by $\Psi_{\vec{n}}$ where $\vec{n}=(a,b,1)$.	If $a=0$ we let $\Psi_{\vec{n}} = (x_0+\zeta,y_0,z_0)$. Hence, to any $(a,b) \in \C^2$ corresponds a unique Legendrian curve (up to a reparametrisation) lying in the hyperplane $\Pi$ with complex normal vector $\vec{n}  =(a,b,1)$. For Legendrian curves of a similar type in $\C^{2n+1}$, $n \in \N$, see \cite[Proposition 6.1]{Alarcon2017}.

	\subsection{Exposing boundary points with hitting}
	
	An \emph{admissible set} (see e.\ g.\ \cite[Definition 1.12.9.]{mincplx} and the discussion therein) in an open Riemann surface $\r$ is a compact subset $S = K \cup \Gamma$ where $K$ is a finite union of pairwise disjoint smoothly bounded closed domains $K_i$ and $\Gamma$ is a finite union of pairwise disjoint smooth Jordan arcs $\gamma_j$. We further require that the interior of every $\gamma_j$ is disjoint from $K$. Thus, if $\gamma_j$ intersects some $K_i$, this may only happen in an endpoint (or both) of $\gamma_j$. If this is the case, we require that the intersection $\gamma_j \cap bK_i$ is transverse. Note that every compact bordered Riemann surface may be viewed as an admissible set in some open Riemann surface.
	
	A map $f\colon S \to X$ of an admissible set into a smooth manifold $X$ is of class $\c^r(S)$ if it can be extended to a map of class $\c^r(U)$ where $U$ is an open neighbourhood of $S$. If $X$ is a complex manifold, such a map belongs to the class $\a^r(S)$ if it is of class $\c^r(S)$ and is holomorphic in the interior $\Int S = \Int K$. A \emph{generalised Legendrian curve of class $\a^r(S)$} in $\C^{2n+1}$ is a map $f\colon S \to \C^{2n+1}$ of class $\a^r(S)$ into $\C^{2n+1}$ that satisfies $f^*\alpha = 0$ along $S$.
	
	A frequently encountered construction is the following. Given an admissible set $S=K\cup \Gamma$ and a generalised Legendrian curve $f\colon K \to \C^{2n+1}$ we wish to extend $f$ to a generalised Legendrian curve $F\colon S \to \C^{2n+1}$. This is always possible since any two points in a connected complex contact manifold may be connected by a Legendrian path. Even more, any path in a complex manifold between two points $p$ and $q$ may be approximated with a Legendrian one. If the given path is already Legendrian near the endpoints, one may also fix the derivatives at the endpoints of the approximating path, see \cite[Theorem A.6]{Alarcon2017}. Hence, it is possible to ensure e.\ g.\ that the image $f'(\gamma_i)$ of every path $\gamma_i \in \Gamma$ intersects $f(K) = f'(K)$ in a single endpoint of the arc $f'(\gamma_i)$.
	
	Another crucial step in the proof is being able to use \cite[Theorem 2.3]{Forstneric2009}. For our purposes the following more precise version is needed, as presented in \cite[Theorem 6.7.1]{mincplx}.
	
	\begin{thm}
		\label{fwbdexpose}
		Suppose $D$ is a relatively compact smoothly bounded domain in an open Riemann surface $\r$. Let $\gamma_1,\ldots,\gamma_k$ be pairwise disjoint smooth paths in $\r$ such that $D' = D \cup \bigcup_{i=1}^k \gamma_i$ forms a connected admissible set in $\r$ and for every $i=1,\ldots,k$ let $c_i$ be a point in the relative interior of the arc $\gamma_i$. For every $i=1,\ldots,k$ denote by $a_i$ the unique endpoint of $\gamma_i$ that intersects $bD$ and denote by $b_i$ the other endpoint of $\gamma_i$. Let $U_i' \Subset U_i$ be neighbourhoods of the point $a_i$ and let $V_i$ be a neighbourhood of the arc $\gamma_i$ for every $i=1\ldots,k$. Let $\Lambda \subset \Int D$ be a finite set and let $m\colon \Lambda \to \N \cup \{0\}$ be a function. There exists a smooth conformal diffeomorphism $\phi\colon \overline{D} \to \phi(\overline{D}) \subset \r$, arbitrarily close to the identity map on $D \backslash \bigcup_{i=1}^kU_i$ such that
		\begin{enumerate}
			\item $\phi|_D \colon D \to \phi(D)$ is biholomorphic,
			
			\item $\phi$ agrees with the identity map to order at least $m(p)$ at every point $p \in \Lambda$, and
			
			\item for every $i=1\ldots,k$ the following holds:
			\[
			\phi(D\cap U_i') \subset U_i \cup V_i, \quad
			\phi(a_i) = b_i, \quad \text{and} \quad
			c_i \in \phi(\Int D).
			\]	
		\end{enumerate}
	\end{thm}

	\subsection{Riemann-Hilbert problem with interpolation}
	
	The following is a more precise statement of \cite[Theorem 3.3]{Alarcon2017} by adding interpolation. The proof stays almost the same, thus we only explain it briefly. Note that $\dist(p,Q)$, where $p \in \C^n$ is a point and $Q \subset \C^n$ is a subset, denotes the standard Euclidean distance from $p$ to $Q$, that is
	\[
	\dist(p,Q) = \inf_{q \in Q} \|p-q\|.
	\]
	
	\begin{thm}
		\label{rhinterpolation}
		Suppose $M$ is a compact bordered Riemann surface, $I \subset bM$ is an arc which is a proper subset of some boundary component of $M$, $f=(\bx,\by,z) \colon M \to \C^{2n+1}$ is a Legendrian curve of class $\a^1(M)$ and for every $u \in bM$ the map
		\[
		\overline{\D} \ni v \mapsto F_u(v) = (X(u,v), Y(u,v), Z(u,v)) \in \C^{2n+1}
		\]
		is a Legendrian disk of class $\a^1(\overline{\D})$, depending continuously on $u \in bM$ and such that $F_u(v) \equiv f(u)$ for every $u \in bM \backslash I$.
		
		For every $\varepsilon > 0$ and every neighbourhood $U \subset M$ of $I$ there exists a holomorphic Legendrian curve $H\colon M \to \C^{2n+1}$ and a neighbourhood $V \Subset U$ of $I$ with a smooth retraction $\rho \colon V \to V \cap bM$ such that the following hold:
		
		\begin{enumerate}[label=\roman*)]
			\item $\sup \{|H(u) - f(u)| : u \in M \backslash V\} < \varepsilon$,
			
			\item $\dist(H(u), F(u, \T)) < \varepsilon$ for every $u \in bM$, and
			
			\item $\dist(H(u), F(\rho(u), \D)) < \varepsilon$ for every $u \in V$.
		\end{enumerate}
		
		Moreover, given a finite set $\Lambda \subset \Int M \backslash U$, the map $H$ may be taken to agree with $f$ to a given finite order $m \in \N$ for every $p \in \Lambda$. 
	\end{thm}
	
	\begin{proof}
		Note that the only novelty with respect to \cite[Theorem 3.3]{Alarcon2017} is the last statement. In the first step, approximate the $(x,y)$-components of $f$ using Mergelyan's theorem with jet-interpolation to ensure the approximants agree with the components of the starting curve to order at least $m$. When constructing the perturbation family for the spray $\tilde{x}$ (see \cite[Equation 3.5]{Alarcon2017}), we also add functions $h_p$ for every $p \in \Lambda$ in order to control the values of the $z$-component at the points $p \in \Lambda$; this is done analogously as in the proof of \cite[Lemma 3.1]{Svetina2024}. When solving the Cousin I problem, we again use \cite[Proposition 5.9.2]{Forstneric2017a}, but now we also use property 3.\ in the proposition to obtain the stated interpolation property.
	\end{proof}
	
	\subsection{Increasing boundary distance}
	\label{cpltMetric}
	In the proof we will use a version of \cite[Lemma 4.1]{Alarcon2015}, adapted to the setting of Legendrian curves. The lemma enables one to approximate a given Legendrian curve, defined on a compact bordered Riemann surface, in uniform topology by a new Legendrian curve such that the distance from a fixed point in the interior of the surface to its boundary is increased. We first state some definitions.
	
	Let $M$ be compact bordered Riemann surface and view it as a smoothly bounded subset of an open Riemann surface $\widetilde{M}$. Once and for all fix a holomorphic immersion $\vartheta \colon \widetilde{M} \to \C$, provided by the theorem of Gunning and Narasimhan (see \cite{Gunning1967}), and suppose $f\colon S \to \C^N$ is a mapping of class $\c^1(S,\C^N)$ on a subset $S \subset X$, that is, $f$ may be extended to a $\c^1$ mapping on an open neighbourhood of $S$. Since $\dd \vartheta$ furnishes a trivialisation for the tangent bundle $T\widetilde{M}\cong \widetilde{M} \times \C$ we may view the differential $\dd f$ as a family of $\R$-linear mappings $\dd f_p \colon  \C \to \C^N$ continuously parametrised by $p \in S$. We define the $\c^1$ norm of $f$ on $S$ by
	\[
	\|f\|_{1,S} := \sup_{p\in S} \max \left\{\|f(p)\|, \|\dd f_p\|\right\},
	\]
	where $\|\cdot\|$ is the standard Euclidean norm on $\R^n$. Meanwhile, the $\c^0$ norm of $f$ on $S$ is given by $\|f\|_{0,S} = \sup_{p\in S}\|f(p)\|$. Note that for a general $S$ both $\|\cdot\|_{0,S}$ and $\|\cdot\|_{1,S}$ may be infinite, hence they are genuine norms only when $S$ is compact.
	
	The function $\dist_f$ on $M \times M$, induced by an $\a^r(M)$-map $f$, for $r \geq 1$, is given by
	\[
	\dist_f(p,q) = \inf_{\gamma \in \Gamma} \int_0^1 \bigg\| \dv{t}f(\gamma(t)) \bigg \| \dd t,
	\]
	where $\Gamma$ is the set of piecewise continuously differentiable paths $\gamma\colon [0,1] \to M$ with endpoints $p,q \in M$. Note that $\dist_f$ defines a metric on $M$ provided $f$ is nonconstant, since the set of critical points of $f$ is discrete. The distance through $f$ from a given point $p \in \Int M$ to the boundary $bM$ is given by
	\[
	\dist_f(p, bM) := \inf_{q \in bM} \dist_f(p,q).
	\]
	If $g\colon \Int M \to \C^n$ is a holomorphic map on the associated bordered Riemann surface we may define $\dist_g$ on $\Int M$ in the same way. An embedding $g$ is called \emph{complete}, provided $\dist_g$ is a complete metric on $\Int M$. If $g$ is obtained as the limit $\lim_{n \to \infty} f_n$ of holomorphic maps $f_n \colon M \to \C^n$ in the compact-open topology, then $g$ is complete if and only if $\dist_{f_n}(p,bM) \to \infty$ as $n \to \infty$ for any given point $p \in \Int M$.
	
	\begin{lm}
		\label{bdincrease}
		Suppose $f \colon M \to \C^{2n+1}$ is a holomorphic Legendrian curve, defined on a compact bordered Riemann surface $M$. Given a finite set $\Lambda \subset \Int M$, a point $p_0 \in \Int M$, a number $m \in \N$, and a couple of numbers $\epsilon, \lambda > 0$, there exists a holomorphic Legendrian curve $F \colon M \to \C^{2n+1}$ satisfying the following:
		\begin{enumerate}[label=\alph*)]
			\item $\|F - f \|_{0,M} < \epsilon$,
			
			\item $\dist_F(p_0, bM) > \lambda$, and
			
			\item $F$ agrees with $f$ to order at least $m$ at every point $p \in \Lambda$.		
		\end{enumerate}
		Moreover, if $\dd f(p) \neq 0$ for every point $p \in \Lambda \cap \{m\geq 1\}$, then $F$ may be chosen an immersion. If also $f|_\Lambda$ is injective, then $f$ may be chosen an embedding.
	\end{lm}
	
	One proves the lemma using an inductive procedure with the help of the following adaptation of \cite[Lemma 6.3]{Alarcon2017}.
	
	\begin{lm}
		\label{bdinclemma}
		Let $M$ be a compact bordered Riemann surface and let $f \colon M \to \C^{2n+1}$, $n \geq 2$, be a Legendrian immersion of class $\a^1(M)$. Suppose we are given a finite set $\Lambda \subset \Int M$, a function $m \colon \Lambda \to \N \cup \{0\}$, a number $\mu > 0$, a smooth map $\mathfrak{Y} \colon bM \to \C^{2n+1}$, and two numbers $\delta, d > 0$ such that
		\begin{enumerate}[label=(\roman*)]
			\item 
			\label{i}
			$\| f -\mathfrak{Y}\|_{0,bM} < \delta$, and
			\item 
			\label{ii}
			$\dist_f(0,bM) > d$.
		\end{enumerate}
		For any $\mu > 0$, the map $f$ may be approximated uniformly on compacts in $\Int M$ by Legendrian curves $\tilde{f} \colon M \to \C^{2n+1}$ of class $\a^1(M)$ satisfying the following properties:
		\begin{enumerate}[label=(\Roman*)]
			\item $\| \tilde{f} - \mathfrak{Y}\|_{0,bM} < \sqrt{\delta^2 + \mu^2}$,
			\item $\dist_{\tilde{f}}(p_0,bM) > d + \mu$ for a given point $p_0 \in \Int M$, and
			\item $\tilde{f}$ agrees with $f$ to order at least $m(p)$ at every point $p \in \Lambda$.
		\end{enumerate}
		If $f$ is an immersion near every point $p \in \Lambda \cap \{m\geq 1\}$, then $\tilde{f}$ may be chosen an immersion. If, furthermore, $f|_\Lambda$ is injective, then $\tilde{f}$ may be made an embedding.
	\end{lm}
	
	\begin{rmk}
		Note that \cite[Lemma 6.3]{Alarcon2017} is a technical result used to establish \cite[Theorem 1.2]{Alarcon2017} and is stated only for the case where $M = \overline{\D}$ is a disk. The generalisation to an arbitrary compact bordered Riemann surface is obtained by using Theorem \ref{rhinterpolation} instead of \cite[Lemma 3.2]{Alarcon2017}.
	\end{rmk}
	
	\begin{proof}[Proof of Lemma \ref{bdinclemma}]
		For simplicity of notation we present the proof in the case $n=1$. The case $n > 1$ is obtained in the same way. The first part, i.\ e.\ \cite[Claim 6.4]{Alarcon2017}, is done in an essentially the same way as in the cited paper by working on each boundary component of $M$ separately.
		
		View $M$ as a smoothly bounded compact domain in an open Riemann surface $\widetilde{M}$. Let $B_1,\ldots, B_K$ be the connected components of the boundary $bM$ in $\widetilde{M}$. Note that each $B_i$ is a compact smooth manifold without boundary of real dimension 1, thus we may identify it with the circle $S^1$. Let $g=(g_1,g_2,g_3) := f - \mathfrak{Y}$ and assume the product $g_1g_3$ has no zeroes on $bM$, which is achieved by a small deformation of $\mathfrak{Y}$ if necessary. Condition \ref{ii} in the lemma implies that there is a number $d_0 > d$ such that $\dist_f(p_0,bM) > d_0$. We then choose a smoothly bounded neighbourhood $U_{bM}$ for the boundary $bM$ in $M$ such that $U_{bM}$ deformation retracts onto $bM$ and $\dist_f(p_0,bU_{bM}) > d_0$ holds (note that $bU_{bM} \subset \Int M$ and $M_0 := M \backslash U_{bM}$ is a compact bordered Riemann surface with boundary $bU_{bM}$).
				
		Choose a number $\epsilon > 0$ to be specified later and let $m \in \N$ be a large enough number so that on every arc $\alpha_{i,j}$ of the form
		\begin{equation}
			\label{mDef}
		\alpha_{i,j} := 
		\left\{
		\exp(2 i \pi t): t \in \left[ \frac{j-1}{m}, 		\frac{j}{m}\right], 
		j \in \Z_{m}
		\right\}
		\end{equation}
		the following estimates are satisfied for arbitrary $u,u' \in \alpha_{i,j}$:
		\begin{equation}
		\label{mProp}
		|f(u)-f(u')| < \epsilon, \quad
		|f(u) - \mathfrak{Y}| < \delta_0, \quad
		|\mathfrak{Y}(u) - \mathfrak{Y}(u')| < \epsilon.	
		\end{equation}
		For every $i,j$ set $u_{i,j} = \exp(2i\pi j/m)$, recall that we have identified each boundary component $B_i$ with the circle $S^1$. As in \cite[Lemma 6.3]{Alarcon2017}, we denote for a given vector $w \in \C^3\backslash \{0\}$ by $\pi_w$ the Hermitian orthogonal projection $\pi_w \colon \C^3 \to \C w \subset \C^3$ onto the complex line $\C w$. We proceed by proving an analogue of \cite[Claim 6.4]{Alarcon2017}. For every $i,j$ choose a point $q_{i,j} \in \alpha_{i,j}$ to be specified later on.
		
		For every $u_{i,j}$ choose a smooth path $\gamma_{i,j}$ in $\widetilde{M}$ which intersects $M$ transversely and only in the point $u_{i,j}$. Label the other endpoint of $\gamma_{i,j}$ by $v_{i,j}$. By a slight deformation of the arcs $\gamma_{i,j}$ if necessary we may assume these arcs are pairwise disjoint. Note that
		\[
		M' = M \cup \bigcup_{i=1}^K \bigcup_{j=1}^{m} \gamma_{i,j}
		\]
		is an admissible set in the open Riemann surface $\widetilde{M}$. Let $\mathtt{I} = \Z_K \times \Z_l$. Choose a number $c > 0$ to be specified later on. We extend $f \colon M \to \C^3$ to a generalised Legendrian curve still denoted by $f \colon M' \to \C^3$ satisfying the properties
		\begin{enumerate}[label=(a\arabic*)]
			\item 
			\label{smalla1}
			$|f(u) - f(u_{i,j})| < c$ for every $u \in \gamma_{i,j}$, 
			\item $(f_1(u)-\mathfrak{Y}_1(u_{i,j}))(f_3(u)-\mathfrak{Y}_3(u_{i,j})) \neq 0$ for every $u \in \gamma_{i,j}$, and
			\item 
			\label{smalla3}
			if $\{J_{a,b}\}_{(a,b)\in \mathtt{I}}$ is any partition of $\gamma_{i,j}$ by Borel measurable subsets, then
			\[
			\sum_{(a,b)\in \mathtt{I}} \mathrm{length}(\pi_{g(q_{a,b})}(f(J_{a,b}))) > 2\mu.
			\]
		\end{enumerate}
		This is achieved by first extending $f$ smoothly over the arcs and then approximating these arcs with Legendrian ones using \cite[Theorem A.6]{Alarcon2017} in the same way as in \cite[Proof of Claim 6.4]{Alarcon2017}. By choosing a small enough $c > 0$ we may furthermore assume the following:
		\begin{enumerate}[label=(a\arabic*)]
			\setcounter{enumi}{3}
			\item 
			\label{smallaa5}
			$|f(u)-f(u')| < \epsilon$ and $|f(u) - \mathfrak{Y}(u')|< \delta_0$ for every pair 
			\[
			(u,u') \in (\gamma_{i,j-1}\cup \alpha_{i,j} \cup \gamma_{i,j}) \times \alpha_{i,j}
			\]
			and every $(i,j) \in \mathtt{I}$.
		\end{enumerate}
		By \cite[Theorem 2.1]{Svetina2024}, there exists a holomorphic Legendrian curve $f' \colon U \to \C^3$, defined on an open neighbourhood $U$ of $M'$ in $\widetilde{M}$ which is arbitrarily close to $f$ in the $\c^1(M')$-topology and such that $f'$ agrees with $f$ to order at least $m(p)$ at every point $p \in \Lambda$; furthermore, $f'$ may be taken an immersion or an embedding, provided $f$ satisfies the corresponding conditions in the statement of the lemma.
		
		Choose open neighbourhoods $U_{i,j}' \Subset U_{i,j} \Subset U\backslash \Lambda$ of the points $u_{i,j}$ and open neighbourhoods $W_{i,j} \Subset U\backslash \Lambda$ of the paths $\gamma_{i,j}$ such that the following hold:
		\begin{enumerate}[label=(b\arabic*)]
			\item 
			\label{smallb1}
			$W_{i,j} \cap M \Subset U_{i,j}'$;
			
			\item $(U_{i,j} \cup W_{i,j}) \cap (U_{i',j'} \cup W_{i',j'}) = \emptyset$ whenever $(i,j) \neq (i',j')$;
			
			\item $|f'(u) - f(u')| < \epsilon$ and $|f'(u) - \mathfrak{Y}(u')| < \delta_0$ for every $u \in W_{i,j-1} \cup U_{i,j-1} \cup \alpha_j \cup U_{i,j} \cup W_{i,j}$, $u' \in \alpha_j$ and every $i,j$;
			
			\item $(f'_1(u)-\mathfrak{Y}_1(u'))(f'_3(u)-\mathfrak{Y}_3(u')) \neq 0$ for every $u \in U_{i,j}$, $u' \in U_{i,j} \cap bM$;
			
			\item 
			\label{smallb5}
			for any path $\gamma \subset U_{i,j} \cup W_{i,j}$ with one endpoint in $U_{i,j}$ and the other endpoint $v_{i,j} \in W_{i,j}$, and any partition $\{J_{a,b}\}_{(a,b) \in \mathtt{I}}$ of $\gamma$ into Borel measurable subsets $J_{a,b}$ the following holds:
			\[
			\sum_{(a,b)\in \mathtt{I}} \mathrm{length}(\pi_{g(q_{a,b})}(f'(J_{a,b}))) > 2\mu.
			\]
			\item 
			\label{smallb6}
			$|f'(u) - f(u')| < c$ for every $u \in U_{i,j} \cup W_{i,j}$ and every $u' \in U_{i,j} \cap bM$.
		\end{enumerate}
		Then, use Theorem \ref{fwbdexpose} to find a conformal diffeomorphism $\phi \colon M \to \phi(M) \Subset U$ such that
		\begin{enumerate}[label=(c\arabic*)]
			\item 
			\label{cbiholo}
			$\phi \colon \Int M \to \phi(\Int M)$ is biholomorphic,
			
			\item 
			\label{cclosetoid}
			$\phi$ is arbitrarily close to the identity map on $L := M \backslash \bigcup_{i,j}U_{i,j}'$,
			
			\item 
			\label{ctangentid}
			$\phi$ is tangent to the identity map to order at least $m(p)$ at every point $p \in \Lambda$,
			
			\item 
			\label{cuijtovij}
			$\phi(u_{i,j}) = v_{i,j}$ for every $i,j$, and
			
			\item 
			\label{cbiguijtovij}
			$\phi(U_{i,j}') \subset U_{i,j} \cup W_{i,j}$ for every $i,j$.
		\end{enumerate}
		If the conditions in the end of the lemma hold and approximation in \ref{cclosetoid} is good enough, the map $F := f' \circ \phi$ may be assumed to be an immersion or an embedding. Moreover, provided the approximations are chosen good enough, it is easy to check that $F$ satisfies the following conditions:
		\begin{enumerate}[label=(A\arabic*)]
			\item $\|F - f \|_{1,L} < \epsilon'$,
			
			\item $|F(u) - f(u')| < \epsilon$ and $|F(u) - \mathfrak{Y}(u)| < \delta_0$ fo every $u,u' \in \alpha_{i,j}$ and every pair $(i,j)$,
			
			\item the function $G\colon u \mapsto (F_1(u)- \mathfrak{Y}_1)(F_3(u)-\mathfrak{Y}_3(u))$ has no zeroes on $bM$, and
			
			\item 
			\label{biga4}
			for any path $\gamma$ with one endpoint in $L$ and the other endpoint in $bM$, and any partition $\{J_{a,b}\}_{(a,b) \in \mathtt{I}}$ of $\gamma$ into Borel measurable sets $J_{a,b}$, the following holds:
			\[
			\sum_{(a,b)\in \mathtt{I}} \mathrm{length}(\pi_{G(q_{a,b})}(F(J_{a,b}))) > 2\mu.
			\]
		\end{enumerate}
		Thus far we have shown that $f$ may be approximated in $\c^0(M)$ topology by a holomorphic Legendrian curve that has increased distance to a neighbourhood of a finite set in the boundary. The goal is now to increase the distance also at the complement of this neighbourhood. This is achieved by using the approximate solution to the Riemann-Hilbert problem, provided by Theorem \ref{rhinterpolation}.
		
		Given $F$ as above we define, following \cite[Lemma 6.3]{Alarcon2017}:
		\[
		\alpha_{i,j}' = \alpha_{i,j} \backslash (U_{i,j-1}\cup U_{i,j}) \subset \alpha_{i,j},
		\]
		and let $D_{i,j} \Subset M \backslash \Lambda$ be a neighbourhood of $\alpha_{i,j}'$ in $M$, admitting a smooth retraction $\rho_{i,j}\colon D_{i,j} \to \alpha_{i,j}'$, such that, given a number $\epsilon'' >0$ to be defined later on, the following conditions hold:
		\begin{enumerate}[label=(B\arabic*)]
			\item $\|F(u) - f(u')\| < \epsilon$ and $\|F(u) - \mathfrak{Y}(u')\| < \delta_0$ for every $u \in D_{i,j}$ and $u' \in \alpha_{i,j}$;
			
			\item $\|F(\rho_{i,j}(u)-F(u))\| < \epsilon''$ for every $u \in D_{i,j}$.
		\end{enumerate}		
		For every $i,j$ define $\Pi_{i,j}$ to be the complex hyperplane with the complex normal vector $G(u_{i,j})$. Use the construction in Section \ref{propEmbHyperplane} to find a continuous map 
		\[
		\Psi_{i,j}\colon (D_{i,j}\cap bM) \times \C \to \C^3
		\]
		such that the following conditions hold for every $u \in D_{i,j}\cap bM$:
		\begin{itemize}
			\item $\Psi_{i,j}(u,0) = F(u)$,
			
			\item $\Psi_{i,j}(u,\cdot) \colon \C \to \C^3$ is a proper Legendrian embedding, and
			
			\item $\Psi_{i,j}(u,\C) \subset F(u) + \Pi_{i,j}$.
		\end{itemize}
	
		We then use the maps $\Psi_{i,j}$ to construct a continous map $H_{i,j}\colon (D_{i,j} \cap bM) \times \overline{\D} \to \C^3$ satisfying
		\begin{enumerate}[label=(C\arabic*)]
			\item $H_{i,j}(u,0) = F(u)$ for every $u \in D_{i,j} \cap bM$,
			
			\item $H_{i,j}(u,\cdot) \colon \overline{\D} \to \C^3$ is a Legendrian disk, holomorphic on a neighbourhood of $\overline{\D}$ in $\C$,
			
			\item $H_{i,j}(u,\overline{\D}) \subset F(u) + \Pi_{i,j}$,
 			
			\item $\|H_{i,j}(u,\zeta) - F(u)\| = \mu$ for every $u \in \alpha_{i,j}'$ and every $\zeta \in b\overline{\D}$,
			
			\item $\|H_{i,j}(u,\zeta) - F(u)\| \leq \mu$ for every $u \in D_{i,j}\cap bM$	and every $\zeta \in \overline{\D}$, and
			
			\item if $u \in D_{i,j}$ is an endpoint of the arc $D_{i,j} \cap bM$, then $H_{i,j}(u,\cdot) \equiv F(u)$.
		\end{enumerate}
		In order to fulfill properties (C3) and (C4) we need to ensure that the map $\Psi_{i,j}(q,\cdot)$ intersects the ball $B(q,\mu) \subset \C^3$ (with centre $q$ and radius $\mu$) transversely for every $q \in D_{i,j} \cap bM$. This is achieved by suitably choosing the parameters $m \in \N$ (see \ref{mDef}), $\epsilon>0$ (see \ref{mProp}), $c > 0$ (see properties \ref{smalla1} -- \ref{smallaa5}), and the neighbourhoods $U_{i,j}'\Subset U_{i,j}$, $W_{i,j}$ (see properties \ref{smallb1} -- \ref{smallb6}) and perhaps slightly perturbing the map $\mathfrak{Y} \colon bM \to \C^{2n+1}$ as follows. 
		
		Fix a $q \in bM$. Since $\Psi_{G(q)}(q,\cdot)$ is nonconstant it follows that there is a $\mu_q' \in (\mu-(d_0-d), \mu]$ such that $\Psi_{G(q)}(q,\cdot)$ intersects the boundary $bB(q,\mu_u')$ transversely, thus the preimage $\Psi_{G(q)}(q,bB(q,\mu_q'))$ is a closed real dimension 1 submanifold in $\C$, i.\ e.\ a disjoint union of circles. We choose the one that bounds a simply connected domain in $\C$ containing the origin.	Since transversality is an open condition, it follows by continuity of $\Psi$ and compactness of $bM$ that there exists a neighbourhood $\Theta_q$ of $f(q)$ in $\C^3$ such that $\Psi_{G(q)}(v,\C)$ intersects $B(v,\mu_u')$ transversely for every $v \in \Theta_q$.
		
		For every boundary component $B_i$, $i=1,\ldots,K$, choose finitely many points $q_{i,j'} \in B_i$, $j'=1,\ldots,m'$ such that the neighbourhoods $\Theta_{i,j'} := \Theta_{q_{i,j'}}$ for $j'=1,\ldots,m'$ cover $f(bB_i)$ for every $i=1,\ldots,K$. When splitting the boundary components $B_i$ into arcs $\alpha_{i,j}$ we choose the number $m \in \N$ big enough and $\epsilon > 0$ small enough such that \ref{mDef} and \ref{mProp} imply that $f(\alpha_{i,j})$ is contained in some $\Theta_{i,j'}$ for every $i = 1,\ldots, K$ and every $j = 1,\ldots,m$. We relabel the points $q_{i,j'}$ to $q_{i,j}$ and the corresponding neighbourhoods $\Theta_{i,j'} \ni f(q_{i,j'})$ to $\Theta_{i,j}$ so that $\alpha_{i,j} \subset \Theta_{i,j}$ holds for every $i=1,\ldots,K$ and every $j = 1,\ldots,m$.	Next, we choose $c>0$ small enough such that the property \ref{smalla1} implies that $f(\gamma_{i,j}) \subset \Theta_{i,j-1} \cap \Theta_{i,j}$ for every pair $(i,j)$.
		
		Note that condition \ref{smallb6} and a choice of a good enough approximation in \ref{cclosetoid}, \ref{cbiguijtovij} imply that
		\[
		F\left(W_{i,j-1} \cup U_{i,j-1} \cup \alpha_j \cup U_{i,j} \cup W_{i,j}\right) \subset \Theta_{i,j}
		\]
		for every choice of indices $i,j$. In particular, $F(\alpha_{i,j}) \subset \Theta_{i,j}$ holds for every pair of indices $i,j$. It follows that the family of equations in (C4) and (C5) parametrises a smooth family of closed simply connected domains in $\C$. It is standard (see e.\ g.\ \cite[Theorem II.5.2]{Goluzin1969}) that this family may be reparametrised to a smooth family of holomorphic maps $\overline{\D} \to \C^3$. After suitably shrinking these disks for parameters near the endpoints of $\alpha_{i,j}'$ we obtain the desired mapping $H_{i,j}\colon D_{i,j}\cap bM \to \C^3$.
				
		We assemble these mappings into a single continous map $H \colon bM \times \overline{\D} \to \C^3$ defined as follows:
		\begin{equation}
			\label{transverseFamily}
		H(u,\zeta) = 
		\begin{cases}
			H_{i,j}(u,\zeta), & \text{if } u \in D_{i,j}\cap bM \text{ for some } i,j;\\
			F(u), & \text{ if } u \in bM \backslash \bigcup_{i,j}D_{i,j}.	
		\end{cases}
		\end{equation}
		Fix a number $\epsilon''' > 0$ to be specified later, let $L' := M \backslash \bigcup_{i,j}D_{i,j}$, and use Theorem \ref{rhinterpolation} to obtain a holomorphic Legendrian curve $\tilde{f}\colon M \to \C^3$ satisfying
		\begin{enumerate}[label=(D\arabic*)]
			\item $\dist(\tilde{f}(u), H(u,bM)) < \epsilon'''$ for every $u \in bM$,
			
			\item $\dist(\tilde{f}(u), H(\rho_{i,j}(u), \overline{\D})) < \epsilon'''$ for every $u \in \bigcup_{i,j}D_{i,j}$,
			
			\item $\|\tilde{f} - F\|_{1, L'} < \epsilon'''$, and
			
			\item $\tilde{f}$ agrees with $F$ (and hence with $f$) to order at least $m(p)$ at every point $p \in \Lambda$.
		\end{enumerate}
		Provided that appropriate conditions in the statement of the lemma hold, $\tilde{f}$ may be taken an immersion or an embbedding. One then chooses $\epsilon, \epsilon', \epsilon''$ and $\epsilon'''$ small enough and shows that $\tilde{f}$ satisfies the conclusion of the lemma in the same way as in \cite[Lemma 6.3]{Alarcon2017}.
	\end{proof}
	
	\begin{proof}[Proof of Lemma \ref{bdincrease}]
		We follow the proof of \cite[Lemma 4.1]{Alarcon2015}. Choose numbers $d_0$ and $\delta_0$ such that $0 < d_0 < \dist_f(p_0,bM)$ and $0 < \delta_0 < \epsilon$. Define
		\[
		c = \frac{\sqrt{6(\epsilon^2-\delta_0^2)}}{\pi},
		\]
		and note that for
		\[
		d_j := d_{j-1} + \frac{c}{j}, \quad \delta_j = \sqrt{\delta_{j-1}+\frac{c^2}{j^2}},
		\]
		it holds that $d_j \to \infty$ and $\delta_j \to \epsilon$ as $j \to \infty$. Set $f_0 = f$ and inductively define a sequence of holomorphic Legendrian curves $f_j \colon M \to \C^{2n+1}$ satisfying the following properties:
		\begin{enumerate}[label=\roman*$_j$)]
			\item $\|f_j - f\|_{0,bM} < \delta_j$,
			
			\item $\dist_{f_j}(p_0, bM) > d_j$, and
			
			\item $f_j$ agrees with $f$ to order at least $m$ at every point $p \in \Lambda$.
		\end{enumerate}
		Suppose such a sequence exists. As explained above the properties i$_j$) guarantee that every $f_j$ in the sequence is $\epsilon$-close to $f$ on $bM$. By the maximum principle this implies $\|f_j-f\|_{0,M}< \epsilon$ for every $j \in \N$. Meanwhile, properties ii$_j$) guarantee that $\dist_{f_j}(p_0,bM) > d_j > \lambda$ for every large enough $j \in \N$. By iii$_j$) every $f_j$ in the sequence agrees with $f$ to a given finite order at the points in $\Lambda$. If the last condition in the lemma holds, we may use \cite[Theorem 3.5]{Svetina2024} to correct a given $f_j$ to an embedding. Thus, the curve $F$ satisfying the conclusion of the lemma, is obtained by choosing a large enough $j \in \N$ and setting $F = f_j$.
		
		The existence of the sequence is on the other hand a simple consequence of Lemma \ref{bdinclemma} by using the following data at step $j-1$: $f=f_{j-1}$, $\mathfrak{Y} = f|_{bM}$, $\delta = \delta_{j-1}$ and $\mu = c/j$. The lemma is now proved.
	\end{proof}
	
	\begin{rmk}
		When constructing the mapping $H$ in the case $n \geq 2$ we could of course take linear Legendrian disks lying in hyperplanes $\Pi_{i,j}$ as described in Section \ref{propEmbHyperplane}. The problem of transversality becomes trivial since the intersection of a (linear) complex line and a ball is always convex, hence simply connected. This is the method used also in the proof of Lemma \ref{cvxInductiveLemma}. There however, the construction breaks down in the case $n=1$ for reasons stated below, see Remark \ref{breakdown}.
	\end{rmk}

	\section{The main Lemma}
	
	The main part of this section is Lemma \ref{cvxInductiveLemma} which is an adaptation to the contact setting of \cite[Lemma 5.2]{Alarcon2015}. The main difference lies in constructing the family of Legendrian disks with the parameter in $bM$. This is accomplished by noting that the real tangent hyperplane to a point in the boundary of the codomain (when the latter is viewed as a smooth manifold) contains a complex hyperplane and using the construction from \cite[Proposition 6.1]{Alarcon2017} explained in Section \ref{propEmbHyperplane}.
	
	We first present a technical result which will be used many times in the constructions below.
	
	\begin{lm}
		\label{inductionStart}
		Suppose $\coD \subset \C^{2n+1}$ is a connected domain in $\C^{2n+1}$ and $f\colon M \to \C^{2n+1}$ is a holomorphic Legendrian curve defined on a compact bordered Riemann surface $M$ such that $f(M) \subset  \coD$. Let $K \subset \Int M$ be a compact set and let $\Lambda' \subset \Int K$ and $\Lambda'' \subset \Int \coD \backslash f(M)$ be finite sets. Suppose $\coL \subset \coD$ is a compact subset such that $f(M\backslash \Int K)\cup \Lambda'' \subset \coD \backslash {\coL}$.
		
		For a given $m \in \N \cup \{0\}$ there exists a Legendrian curve $f'\colon M \to \coD$ arbitrarily close to $f$ in the $\c^1(K)$ topology such that
		\begin{enumerate}[label=\Roman*)]
			\item $f'(M \backslash \Int K) \subset \coD \backslash {\coL}$,
			\item $f'$ agrees with $f$ to order $m$ at every point $p \in \Lambda'$, and
			\item $\Lambda'' \subset f'(\Int M)$.
		\end{enumerate}
	\end{lm}
	
	\begin{proof}
		We may assume $M$ is a smoothly bounded compact domain in an open Riemann surface $\widetilde{M}$. For every $p \in \Lambda''$ choose a smooth Jordan arc $\gamma_p$ in $\widetilde{M}'$ with one endpoint $u_p$ in the boundary $bM$ and the other endpoint $v_p$ in $\widetilde{M}\backslash M$. Moreover, ensure that $\gamma_p$ intersects $M$ in the single point $u_p \in bM$ such that the intersection $u_p \cap bM$ is transverse and that the arcs $\gamma_p$ are pairwise disjoint. Note that $M' = M \cup \bigcup_{p\in \Lambda''} \gamma_p$ is an admissible set in $\widetilde{M}$ and extend $f$ to a generalised Legendrian curve still denoted by $f \colon M' \to \C^{2n+1}$ such that $f(v_p') = p$ for some $v_p' \in \Int \gamma_p$ and $f(\gamma_p) \subset \coD \backslash {\coL}$ for every $p \in \Lambda''$. Moreover, we may assume $f$ is holomorphic near every point $v_p$. Use \cite[Theorem 2.1]{Svetina2024} to approximate $f$ with a holomorphic Legendrian curve $g \colon U \to \C^{2n+1}$, defined on a neighbourhood $U$ of $M'$ in $\widetilde{M}$ such that $g$ agrees with $f$ to order $m$ at every point $p \in \Lambda'$ and such that $g(v_p) = f(v_p) = p$. By shrinking $U$ if necessary we may assume $g(bU)$ does not intersect $\coL$. Now use Theorem \ref{fwbdexpose} to find a conformal diffeomorphism $\phi \colon M \to \overline{U}'$ where $U' \subset U$ is a possibly smaller open neighbourhood of $M'$ in $\widetilde{M}$ such that $\Lambda'' \subset \phi(\Int M)$ and $\phi$ is tangent to the identity map to order $m$ at every point $p \in \Lambda'$. Moreover, $\phi$ may be chosen arbitrarily close to the identity map on the compact set $K$. Then, $f' := g \circ \phi$ satisfies the conclusion of the lemma.
	\end{proof}
	
	Now suppose $\coD \Subset \R^{n+1}$ is a relatively compact smoothly bounded strongly convex domain. It is standard that its boundary $b\coD$ is diffeomorphic to the $n$-sphere. Define $N_\coD\colon b\coD \to \R^{n+1}$ to be the inward pointing unit normal vector field to $b\coD$. Denote by $\kappa_{\coD}^\mmax(p)$ and $\kappa_{\coD}^\mmin(p)$ the maximal and minimal principal curvatures of $b\coD$ at the point $p \in b\coD$ and let
	\[
	\kappa_{\coD}^\mmax = \max_{p\in b\coD} \kappa_{\coD}^\mmax(p), \quad 
	\kappa_{\coD}^\mmin = \min_{p\in b\coD} \kappa_{\coD}^\mmin(p).
	\]
	By pushing $b\coD$ in the direction of $N_\coD$ for time $t < 1/\kappa_{\coD}^\mmax$ we obtain another bounded strongly convex domain $\coD_t$ whose boundary is given by
	\begin{equation}
		\label{mt}
		b\coD_t = \{p+tN_\coD(p):p \in b\coD\}.
	\end{equation}
	It is elementary that the domain $\coD_t$ is well defined and strongly convex for every $t \in (-\infty,\tau)$ for some $\tau > 0$. The minimal and maximal principal curvatures $\kappa_{\coD_t}^\mmin$ and $\kappa_{\coD_t}^\mmax$ are given by equations (see e.\ g.\ \cite[Theorem 1.18]{Bellettini2013})
	\[
	\frac{1}{\kappa_{\coD_t}^\mmin} = \frac{1}{\kappa_{\coD}^\mmin}	- t, \quad
	\frac{1}{\kappa_{\coD_t}^\mmax} = \frac{1}{\kappa_{\coD}^\mmax} - t.
	\]
	It follows immediately that $\tau < 1/\kappa_{\coD}^\mmax$. That in fact $\tau = 1/\kappa_{\coD}^\mmax$ follows by an elementary computation.
	
	\begin{lm}
		\label{cvxInductiveLemma}
		Let $\coL \Subset \coD$ be relatively compact smoothly bounded strongly convex domains in $\C^{2n+1}$ such that $\coD \subset \coL_{-\eta}$ for some $\eta>0$ and let $f\colon M \to \coD$ be a holomorphic Legendrian curve defined on a compact bordered Riemann surface $M$. Assume for a compact set $K \subset \Int M$ we have
		$$
		f(M \backslash \Int K) \subset \coD \backslash \overline{\coL}.
		$$
		Given a number $0 < \delta < \tau_\coD$, a finite set of points $\Lambda' \subset \Int M$, and a number $m \in \N\cup \{0\}$ the curve $f$ can be approximated as closely as desired in $\c^1(K)$-topology by a holomorphic Legendrian curve $F\colon M \to \coD$, satisfying the following properties:
		\begin{enumerate}[label={\alph*)}]
			\item 
			\label{sqrtrazdalja}
			$\|F - f\|_{0,M} < \sqrt{2\eta^2 + 2\eta/\kappa_\coL^{\min}}$,
			
			\item
			\label{Kblizu}
			$\|F-f\|_{0,K} < \varepsilon$,
			
			\item 
			\label{izvenkompakta}
			$F(M\backslash \Int K) \subset \coD \backslash \overline{\coL}$,
			
			\item 
			\label{slikaroba}
			$F(bM) \subset \coD \backslash \overline{\coD }_\delta$, and
			
			\item 
			\label{interiorInterpolation}
			$F$ agrees with $f$ to order at least $m$ at every point $p \in \Lambda'$.
		\end{enumerate}
	\end{lm}
	
	\begin{proof}[Proof of Lemma \ref{cvxInductiveLemma}]
		We first treat the case with no hitting, that is, $\Lambda'' = \emptyset$, and prove a) -- d). Assume $M$ is a smoothly bounded compact domain in an open Riemann surface $\widetilde{M}$ and that $\delta > 0$ is small enough, so that $\overline{\coL} \subset \coD_\delta$. Choose a constant $\varsigma > 0$ such that
		$$
		\overline{\coL}_{-\varsigma} \subset \coD_\delta, \quad f(M \backslash \Int K) \subset \coD \backslash \overline{\coL}_{-\varsigma}.
		$$
		Pick a constant $c > 0$, which will be determined later, and for every $x \in b\coL$ set
		\begin{equation}
			\label{cdefined}
			B_x = b\coL \cap \mathbb{B}_x(c),
		\end{equation}
		where $\mathbb{B}_x(c)$ is an open Euclidean ball in $\C^{2n+1}$ with centre $x$ and radius $c$. Set
		$$
		O_x = \coD \cap \{y - tN_\coL(y)\, : \, y \in B_x, \, t > \varsigma\} \subset \coD \backslash \overline{\coL}_{-\varsigma}.
		$$
		Assume $c>0$ is small enough, so that $B_x$ is a topological open ball and
		\begin{equation}
			O_x \subset \widetilde{O}_x = \{y \in \coD\, : \,
			\langle y-x, \, N_\coL(x)\rangle < - \varsigma/2\}
			\subset \coD \backslash \overline{\coL}_{-\varsigma/2}
			\quad \text{for every} \; x \in b\coL.
		\end{equation}
	
	\begin{figure}[h]
		\begin{tikzpicture}
			\node[anchor=south west,inner sep=0] (image) at (0,0) {\includegraphics[height=5.2cm]{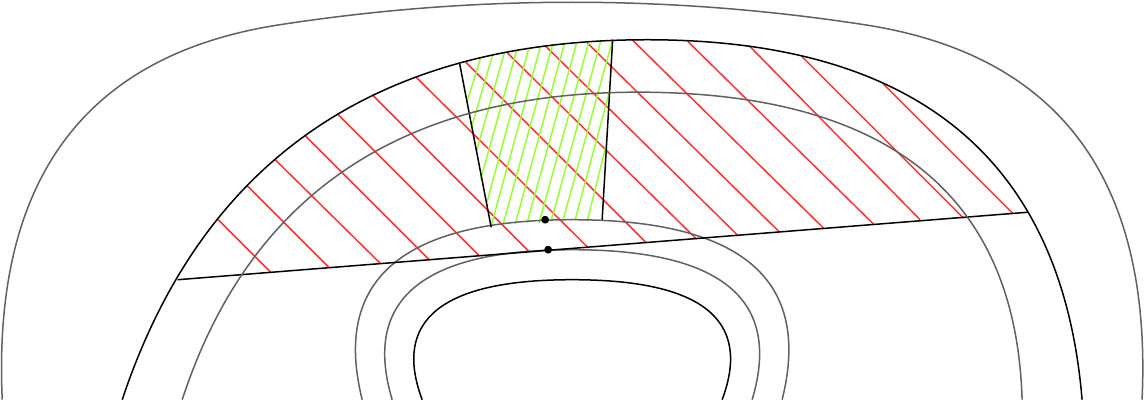}};
			\begin{scope}[x={(image.south east)},y={(image.north west)}]
				\node at (0.04,0.05) {$\coL_{-\eta}$};
				\node at (0.13,0.057) {$\coD$};
				\node at (0.183,0.05) {$\coD_{\delta}$};
				\node at (0.286,0.05) {$\coL_{-\varsigma}$};
				\node at (0.388,0.057) {$\coL$};
				\node [fill=white, rounded corners] at (0.48,0.65) {$O_x$};
				\node [fill=white, rounded corners] at (0.57,0.31) {$x-\varsigma N_\coL(x)/2$};
				\node [fill=white, rounded corners] at (0.7,0.55) {$\widetilde{O}_x$};
			\end{scope}
		\end{tikzpicture}
		
		\caption{The sets $\coD$, $\coL$, $O_{x}$, and $\widetilde{O}_x$.}
	\end{figure}

		Denote by $\alpha_1, \ldots, \alpha_k$ the (connected) Jordan curves forming the boundary components of $M$.
		Since the sets $O_x$ form an open cover for $f(bM)$ we may find a finite subcover consisting of sets $O_{i,j} := O_{x_{i,j}}$, where $i=1,\ldots,k$ and $j \in \Z_l = \Z / l\Z$ for some $l \in \N$.


		Choose smooth Jordan arcs $\alpha_{i,j} \subset bM$ such that 
		\begin{itemize}
			\item $\alpha_i = \cup_{j=1}^l \alpha_{i,j}$ for every $i = 1,\ldots,k$;
			
			\item $\alpha_{i,j} \cap \alpha_{i,j+1} = \{p_{i,j}\}$ for every $i$ and $j \in \Z_l$, otherwise the arcs are pairwise disjoint,
			
			\item $f(\alpha_{i,j}) \subset O_{i,j}$ for some $x_{i,j} \in b\coL$ and $\nu \in 0,\ldots ,N$.
		\end{itemize}
		
		We now expose the boundary points $p_{i,j}$ as follows. For each such point choose a smooth Jordan arc $\gamma_{i,j} \subset \widetilde{M}$ with one endpoint $p_{i,j} \in bM$ and the other endpoint $q_{i,j} \in \widetilde{M}\backslash M$ that intersects $M$ transversely in $p_{i,j}$ and does not intersect it otherwise. Then, perturb the arcs $\gamma_{i,j}$ if necessary such that the family of arcs is pairwise disjoint. Note that
		\[
		M' = M \cup \bigcup_{i,j} \gamma_{i,j}
		\subset \widetilde{M}
		\]
		forms an admissible set in $\widetilde{M}$. By reparametrising if necessary we may assume $\gamma_{i,j}(0) = p_{i,j}$ and $\gamma_{i,j}(1) = q_{i,j} \in \widetilde{M}\backslash M$. Extend $f$ to a generalised Legendrian curve $f' \colon M' \to \coD$ such that
		
		\begin{enumerate}[label=\roman*)]
			\item
			\label{fextend}
			$f'|_M = f$,
			
			\item 
			\label{gammaintersect}
			$f'(\gamma_{i,j}) \subset O_{i,j} \cap O_{i,j+1}$, and
			
			\item
			\label{pushboundary}
			$f'(q_{i,j}) \subset \coD \backslash \overline{\coD}_{\delta/2}$ for all $i,j$.
		\end{enumerate}

		Approximate $f'$ with a holomorphic Legendrian curve $g\colon U \to \C^{2n+1}$ defined on a neighbourhood $U$ of $M'$ in $\widetilde{M}$ which we shrink if necessary to ensure $g(U) \subset \coD$ and $g(bU) \subset \coD \backslash \overline{\coL}_{-\varsigma}$. Moreover, we may assume $g$ agrees with $f$ to order at least $m$ at every point $p \in \Lambda'$. Property i) ensures $g$ approximates $f$ on $M$. If this approximation is good enough, by properties ii) and iii) there exist open neighbourhoods $U_{i,j}$ of the arcs $\alpha_{i,j}$, neighbourhoods $V_{i,j}$ of the arcs $\gamma_{i,j}$ and neighbourhoods $W_{i,j}' \Subset W_{i,j} \Subset U_{i,j} \cap U_{i,j+1}$ of the points $p_{i,j}$, all subsets of $U\backslash K$, such that
		
		\begin{enumerate}[label=\roman*)]
			\setcounter{enumi}{3}
			
			\item $V_{i,j} \cap M \Subset W'_{i,j}$,
			
			\item $g(V_{i,j} \cup U_{i,j} \cup V_{i,j+1}) \subset O_{i,j}$,
			
			\item \label{tri} $g(q_{i,j}) \in \coD \backslash \overline{\coD}_{\delta/2}$.
		\end{enumerate}
		
		By Theorem \ref{fwbdexpose} there exists a conformal diffeomorphism $\phi\colon \overline{M} \to \phi(\overline{M})$, biholomorphic on $\Int M$, arbitrarily close to the identity map on $M\backslash \bigcup_{i,j}W_{i,j}'$ satisfying
		\begin{enumerate}[label=\Roman*)]		
			\item $\phi(\overline{M}) \subset U$,
			
			\item $\phi(W_{i,j}') \subset W_{i,j} \cup V_{i,j}$,
			
			\item $\phi(p_{i,j}) = q_{i,j}$ for every $(i,j) \in \Z_k \times \Z_l$,
			
			\item $\phi$ agrees with the identity map to order at least $m$ at every point $p \in \Lambda'$, and
			
			\item $\gamma_{i,j} \subset \phi(\overline{M})$ for every $i,j$.
		\end{enumerate}
		By I) the composition $G := g \circ \phi$ is well defined. Note that $G$ approximates $f$ on $M\backslash \bigcup_{i,j}W_{i,j}'$ while satisfying $G(M\backslash \Int K) \subset \coD \backslash \overline{\coL}_{-\varsigma}$. Since $g$ agrees with $f$ to order at least $m$ at every point $p \in \Lambda'$, property IV) implies the same for $G$. By properties v), \ref{tri}, II) and III) there exists proper connected subarcs $\beta_{i,j} \subset \alpha_{i,j}$ such that 
		\begin{equation}
			\label{Gboundary}
			G(\alpha_{i,j}\backslash \Int \beta_{i,j}) \subset 
			(\coD \backslash \overline{\coD}_{\delta/2}) \cap O_{i,j}.
		\end{equation}
		
		We now use the Riemann-Hilbert method to push $G(\beta_{i,j})$ into $\coD \backslash \overline{\coD}_{\delta/2})$.
		
		\begin{claim}
			\label{claim1}
			There exists a continuous function $H\colon \alpha_{i,j} \times \overline{\D} \to \C^{2n+1}$ such that for every $u \in \alpha_{i,j}$ the map $H(u,\cdot)\colon \overline{\D} \to \C^{2n+1}$ is a holomorphic Legendrian embedding satisfying the properties
			\begin{enumerate}[label=\Alph*)]
				\item $H(u,0) = G(u)$,
				\item 
				\label{bigB}
				$H(u,\overline{\D}) \subset \coD \backslash \overline{\coL}_{-\varsigma/2}$, and
				\item 
				\label{Hboundarydef}
				$H(u,b\overline{\D}) \subset \coD \backslash \overline{\coD}_\delta$.
			\end{enumerate}
		\end{claim}
		
		\begin{proof}[Proof of Claim]
			Recall that $B_{i,j} := B_{x_{i,j}}$ denotes the projection of the set $O_{i,j}$ onto the boundary $b\coL$ and denote by $B_{i,j}'$ the projection of $O_{i,j}$ onto the boundary $b\coL_{-\varsigma}$. Let $\Sigma = T_{x_{i,j}}b\coL$ be the real tangent hyperplane to $b\coL$ at $x_{i,j}$ and note that after shrinking the radius $c>0$ of $B_{i,j}$ if necessary we may assume that $\overline{\coL}_{-\varsigma/2}$ is contained in the open halfspace defined by the hyperplane $x + \Sigma$ for every $x \in O_{i,j}$ since $\overline{\coL}_{-\varsigma/2}$ is strictly convex. Let $\Pi$ be the unique complex hyperplane contained in $\Sigma$ and write
			\[
			\Pi = 
			\left\{
			(\bx,\by,z) \in \C^{2n+1} : \langle (\mathbf{a},\mathbf{b},c) , (\bx,\by,z) - x_{i,j}\rangle = 0
			\right\}
			\]
			where $\langle\cdot,\cdot \rangle$ denotes the standard Hermitian product on $\C^{2n+1}$ and $(\mathbf{a},\mathbf{b},c) \in \C^{2n+1}$ is the complex normal vector to $\Pi$. Denote by $\Pi(u) \subset \C^{2n+1}$ the parallel complex hyperplane to $\Pi$ that intersects the point $G(u) \in G(\alpha_{i,j})$.

			By a slight perturbation we may assume the product $a_1\cdots a_n$ for $\mathbf{a}=(a_1,\ldots,a_n)$ is nonzero while $x + \Pi$ still does not intersect $\overline{\coL}_{-\varsigma/2}$ for every $x \in O_{i,j}$. Choose a number $0 < \delta' < \delta/2$ such that $\beta_{i,j} \subset \coD_{\delta'}$.
			
			For every pair $(a_1,a_2) \in (\C^*)^2$ let $\chi\colon \C \to \C^{2n+1}$ denote the mapping defined by $\zeta \mapsto (a_2\zeta,-a_1\zeta,0,\ldots,0)$ and for every $u \in \beta_{i,j}$ let $\chi_u\colon \C \to \C^{2n+1}$ be the mapping
			\[
			\chi_u\colon \zeta \mapsto G(u)+\chi(\zeta).
			\]
			Note that for every $u \in \Int (\beta_{i,j}\cap \coD_{\delta'})$ the map $\chi_u$ is a linear holomorphic Legendrian embedding of $\C$ into $\C^{2n+1}$ whose image is entirely contained in the hyperplane $\Pi(u)$. In particular, since $G(u)$ is an interior point of $\coD_{\delta'}$ and the latter is convex, the intersection $\chi_u(\C) \cap b\coD_{\delta'}$ is transverse for every $u \in \beta_{i,j}$. It follows that the connected component of $D_u:=\chi_u^{-1}(\coD_{\delta'})$ containing the origin is biholomorphic to the unit disk $\D$. The desired mapping $H$ is obtained by reparametrising the family $\chi_u$.
			
		\end{proof}
		\begin{rmk}
			\label{breakdown}
			Note the similarity of the above construction of the map $H$ with that of a similar mapping \ref{transverseFamily} in the proof of Lemma \ref{bdinclemma}. The crucial distinction lies in the assumption \ref{bigB} in Claim \ref{claim1}. It is for this reason that here we are unable to repeat the method from Lemma \ref{bdinclemma} as we need to control the global position of disks $\chi_u$ which prevents us from splitting the boundary $bM$ into smaller pieces in which the constructed disks are transverse.
		\end{rmk}
		
		By \cite[Theorem 3.3]{Alarcon2017} there exist pairwise disjoint open neighbourhoods $\Omega_{i,j} \subset M\backslash (K \cup \Lambda' \cup \Lambda'')$ of the arcs $\beta_{i,j}$, smooth retractions $\rho_{i,j}\colon \Omega_{i,j} \to \beta_{i,j}$ with $\Omega = \cup_{i,j}\Omega_{i,j}$ and a holomorphic Legendrian curve $F\colon M \to \coD$ such that
		\begin{enumerate}
			\item 
			\label{fapprox}
			$\sup\{|F(u)-f(u)|\} < \varepsilon$ on $M \backslash \Omega$,
			
			\item 
			\label{Hboundary}
			$\dist(F(u), H(u,b\D)) < \varepsilon$ for every $u \in bM$,
			
			\item 
			\label{Htwist}
			$\dist(F(u), H(\rho(u), \D)) < \varepsilon$ for every $u \in \Omega$, and
			
			\item
			\label{Finterpolation}
			$F$ agrees with $f$ to order at least $m$ at every point $p \in \Lambda'$.
		\end{enumerate}
		
		We check that $F$ satisfies the required properties in the same way as in \cite[Lemma 5.2]{Alarcon2015}. Namely, let $p \in bM$. If $p \in \alpha_{i,j} \backslash \Int \beta_{i,j}$ for some $i,j$, then $p \in \coD \backslash \overline{\coD}_\delta$ by \ref{Gboundary} and \ref{fapprox}, provided that $\varepsilon$ is chosen small enough. If, on the other hand, $p \in \beta_{i,j}$ for some $i,j$, then the same holds by \ref{Hboundarydef}, \ref{Hboundary} and the definition of $H$. This proves property \ref{slikaroba}. Next, by shrinking $\Omega_{i,j}$ if necessary to ensure $f(\overline{\Omega_{i,j}}) \subset O_{i,j}$ for all $i,j$,  property \ref{izvenkompakta} holds by the assumption on $f$ and \ref{fapprox} provided that $\varepsilon$ is chosen small enough. It remains to check property \ref{sqrtrazdalja}. If $p \in M \backslash \Omega$ this is trivial if $\varepsilon$ is chosen such that $\varepsilon < \sqrt{2\eta^2 + 2\eta/\kappa_\coL^\mmin}$. Suppose now $p \in \Omega_{i,j}$ for some $i,j$ and note that $1/\kappa_{\coL_{-\eta}}^\mathrm{min} = 1/\kappa_\coL^\mathrm{min}+\eta$, thus $\coL_{-\eta}$ is contained in the ball with center $x_{i,j}+(1/\kappa_\coL^\mathrm{min})N_\coL(x_{i,j})$ of radius $1/\kappa_\coL^\mathrm{min}+\eta$. From $F(p) \in \widetilde{O}_{x_{i,j}}$ we get $\langle F(p)-x_{i,j}, N_\coL(x_{i,j}) \rangle < - \varsigma/2$ and from 
		\[
		F(p) \in B(x_{i,j}+(1/\kappa_\coL^\mmin)N_\coL(x_{i,j}), \eta + 1/\kappa_\coL^\mmin)
		\]
		we get $\|F(p) - (x_{i,j} + 1/\kappa_\coL^\mmin N_\coL(x_{i,j}))\| < \eta + 1/\kappa_\coL^\mmin$. Now compute for $t \in (\varsigma/2, \eta)$:
		\begin{multline}
				\|F(p) - (x_{i,j}-t N_\coL(x_{i,j}))\|^2  = \\
				\shoveright{=\|F(p) - (x_{i,j} + (1/\kappa_\coL^\mmin) N_\coL(x_{i,j})) + 
				(1/\kappa_\coL^\mmin)N_\coL(x_{i,j}) + tN_\coL(x_{i,j}) \|^2=} \\
				\shoveright{=\|F(p) - (x_{i,j} + 1/\kappa_\coL^\mmin N_\coL(x_{i,j}))\|^2}\\
				\shoveright{+2 \langle F(p)-(x_{i,j}+(1/\kappa_\coL^\mmin)N_\coL(x_{i,j})), 
				(1/\kappa_\coL^\mmin +t)N_\coL(x_{i,j})\rangle
				+ \| (1/\kappa_\coL^\mmin+t)N_\coL(x_{i,j})\|^2}\\
				\shoveright{< (\eta+1/\kappa_\coL^\mmin)^2 + 2\langle F(p)-x_{i,j}, (1/\kappa_\coL^\mmin +t)N_\coL(x_{i,j})\rangle -}\\
				\shoveright{2\langle (1/\kappa_\coL^\mmin)N_\coL(x_{i,j}), (1/\kappa_\coL^\mmin +t)N_\coL(x_{i,j})  \rangle
				+ (1/\kappa_\coL^\mmin +t)^2 }\\
				\shoveright{< (\eta+1/\kappa_\coL^\mmin)^2 - 
				\varsigma (1/\kappa_\coL^\mmin + t) - 2/\kappa_\coL^\mmin(1/\kappa_\coL^\mmin+t) + 
				 (1/\kappa_\coL^\mmin +t)^2 = }\\
				= \eta^2 + 2\eta/\kappa_\coL^\mmin - \varsigma/\kappa_\coL^\mmin-t\varsigma + t^2.
		\end{multline}
		Since this inequality holds for all $t \in (\varsigma/2,\eta)$ we note that $t \mapsto t^2-\varsigma t$ is increasing on $(\varsigma/2, \eta)$ to plug in $t = \eta$ on the right hand side, rearrange and take the square root to obtain
		\begin{equation}
			\label{tildeF}
			\|F(p) - (x_{i,j}-t_p N_\coL(x_{i,j}))\| < \sqrt{2\eta^2 + \left(\frac{2}{\kappa_\coL^\mmin}-\varsigma\right)\eta - \frac{\varsigma}{\kappa_\coL^\mmin}} <
			\sqrt{2\eta^2+ \frac{2\eta}{\kappa_\coL^\mmin}} - \epsilon
		\end{equation}
		for some $\epsilon > 0$. Now denote $x=x_{i,j}$, $x' = x_{i,j} + N_\coL(x_{i,j})/\kappa_\coL^\mmax$ and let $C$ be the circle with center $x'$ and radius $1/\kappa_\coL^\mmax$ lying in the real affine plane containing the line segments $L = [x'x]$ and $L'=[x'f(p)]$. Note that $C \subset \overline{\coL}$. Let $y$ be the point of intersection between $C$ and $L'$. Note that the length of $[x'y]$ is strictly less than $c > 0$. Denote by $\alpha$ the angle $\alpha = \angle xx'y$. The triangle $\Delta xx'y$ is isosceles with $|\overline{xx'}| = |\overline{x'y}|= 1/\kappa_\coL^\mmax$ and $|xy| < c$, thus
		\[
		\sin (\alpha/2) < \frac{c \kappa_\coL^\mmax}{2}.
		\]
		It follows that
		\begin{align}
			\begin{split}
				\label{smallf}
				\|f(p) - (x_{i,j}-t_pN_\coL(x_{i,j}))\| &= 
				(1/\kappa_\coL^\mmax +t_p) \tan \alpha \leq\\
				&\leq (1/\kappa_\coL^\mmax + \eta) 
				\frac{2\sin (\alpha/2) \cos (\alpha/2)}{1-2\sin^2(\alpha/2))} \leq \\
				&\leq (1/\kappa_\coL^\mmax + \eta) \frac{c \kappa_\coL^\mmax}{1-c^2(\kappa_\coL^\mmax)^2/4}  =\\
				&= c\cdot \frac{(1+\eta\kappa_\coL^\mmax)}{1-c^2(\kappa_\coL^\mmax)^2/4} =\\
				&= c(1+\eta\kappa_\coL^\mmax) + o(c).
			\end{split}
		\end{align}
		By adding up the terms in \ref{tildeF} and \ref{smallf} and choosing $c$ small enough we get
		\[
		\|F(p) - f(p)\| < \sqrt{2\eta^2+ \frac{2\eta}{\kappa_\coL^\mmin}}
		\]
		proving property \ref{sqrtrazdalja}.
	\end{proof}
	
	\begin{cor}
		\label{bdInterpolation}
		With everything as in the statement of Lemma \ref{cvxInductiveLemma}, suppose we are additionally given a finite set of points $\Lambda'' \subset \coD \backslash \coD_\delta$. Then, there exists a holomorphic Legendrian curve $F\colon M \to \coD$, arbitrarily close to $f$ in the $\c^1(K)$-topology, such that
		\begin{enumerate}[label={\Alph*)}]
			\item 
			\label{izvenkompakta2}
			$F(M\backslash \Int K) \subset \coD \backslash \overline{\coL}$,
			
			\item 
			\label{slikaroba2}
			$F(bM) \subset \coD \backslash \coD_\delta$,
			
			\item 
			\label{interiorInterpolation2}
			$F$ agrees with $f$ to order at least $m$ at every point $p \in \Lambda'$, and
			
			\item 
			\label{hitting2}
			$\Lambda'' \subset F(\Int M)$.
		\end{enumerate}
	\end{cor}
	
	\begin{proof}
		Use Lemma \ref{inductionStart} to approximate $f$ with a holomorphic Legendrian curve $f'\colon M \to \coD$ satisfying \ref{izvenkompakta2}, \ref{interiorInterpolation2}, and \ref{hitting2}. Then, use Lemma \ref{cvxInductiveLemma} to approximate $f'$ with a holomorphic Legendrian curve $F\colon M \to \coD$ satisfying also \ref{slikaroba2}. Note that $F$ approximates $f'$ uniformly on $M$ with the bound given by \ref{sqrtrazdalja} in Lemma \ref{cvxInductiveLemma}, but the same bound obviously cannot hold for the distance of $F$ to the original map $f$.
	\end{proof}
	
	\section{Proofs of main results}
	\label{pfsMain}
	In this section we prove general versions of Theorems \ref{properCvxBasic} and \ref{weakCvx} from the introduction.
	
	We first present Theorem \ref{properCvx} on approximation of holomorphic Legendrian curves in strongly convex domains by proper and complete ones. The proof is by induction using Lemma \ref{cvxInductiveLemma} and the modification of \cite[Lemma 6.3]{Alarcon2017}, namely Lemma \ref{bdinclemma}. At each step of the induction one first uses Lemma \ref{cvxInductiveLemma} to push the boundary of the image of the bordered Riemann surface closer to the boundary of the given convex set while hitting some prescribed points. Then, using Lemma \ref{bdinclemma}, one perturbs the Legendrian curve in order to increase the induced distance on the bordered Riemann surface while fixing a finite set of points. In the limit one obtains a continuous proper curve whose restriction to the interior is a complete Legendrian embedding, hitting a specified discrete set of points in a given convex domain. 
	
	Note that the theorem also provides a uniform estimate on the norm of the difference between the starting curve and the approximant as a function of both the curvature of the boundary of the ambient strongly convex domain and the distance from this boundary of the image of the starting curve. This estimate is used in Theorem \ref{convExhaustionThm} in order to construct proper Legendrian curves in general convex domains.
	
	\begin{thm}
		\label{properCvx}
		Suppose $\coL \Subset \coD$ are relatively compact smoothly bounded strongly convex domains in $\C^{2n+1}$ satisfying $\overline{\coD} \subset \coL_{-\eta}$ for some $\eta > 0$, and suppose $f\colon M \to \C^{2n+1}$ is a Legendrian curve defined on a compact bordered Riemann surface $M$ such that $f(M) \subset \coD$. Let $K \subset \Int M$ be a compact set such that $f(M \backslash \Int K) \subset \coD \backslash \overline{\coL}$, and let $\Lambda \subset \Int K$ be a finite set. Given a number $\varepsilon > 0$ and a number $m \in \N \cup \{0\}$ there exists a continous map $F\colon M \to \overline{\coD}$ such that
		\begin{enumerate}[label=\Roman*)]
			\item $\|F-f\|_{0,K} < \varepsilon$,
			
			\item $\|F-f\|_{0,M} < \sqrt{2\eta^2 + 2\eta/\kappa_{\coL}^\mmin}$,
			
			\item $F(M \backslash \Int K)\subset \overline{\coD}\backslash \overline{\coL}$ and $F(bM) \subset b\coD$,
			
			\item $F|_{\Int M} \colon \Int M \to \coD$ is a proper and complete Legendrian curve, and		
			
			\item $F$ agrees with $f$ to order at least $m$ at every point $p \in \Lambda$.
		\end{enumerate}
		Moreover, if either $m = 0$ or $m \geq 1$ and $\dd f(p) \neq 0$ for every point $p \in \Lambda$, then $F|_{\Int M}$ may be chosen an immersion. If, furthermore, $f|_{\Lambda}$ is injective, then $F|_{\Int M}$ may be chosen an embedding.
	\end{thm}
	
	\begin{proof}
		Let $\epsilon > 0$ be such that $\overline{\coL}_{-\epsilon} \subset \coD$ and $f(bM) \subset \coD \backslash \overline{\coL}_{-\epsilon}$. Since $1/\kappa_{\coL_{-\epsilon}}^\mmin = 1/\kappa_\coL^\mmin + \epsilon$ we obtain the following by a short computation:
		\[
		\sqrt{2(\eta-\epsilon)^2 + 2\frac{\eta-\epsilon}{\kappa_{\coL_{-\epsilon}}^\mmin}} <
		\sqrt{2\eta^2+ \frac{2\eta}{\kappa_\coL^\mmin}}.
		\]
		Choose a strictly decreasing sequence of numbers $\delta_j>0$ such that $f(bM) \subset \coD_{\delta_0}$ and the following holds:
		\begin{equation}
			\label{deltainequality}
			\sqrt{2(\eta-\epsilon)^2 + 2\frac{\eta-\epsilon}{\kappa_{\coL_{-\epsilon}}^\mmin}} +
			\sum_{j=1}^\infty \sqrt{2\delta_{j}^2+2\frac{\delta_{j}}{\kappa_{\coD}^\mmin}} <
			\sqrt{2\eta^2+ \frac{2\eta}{\kappa_\coL^\mmin}}.
		\end{equation}
		
		Let $\coD_j := \coD_{\delta_j}$, fix a point $p_0 \in \Int M$ and let $\delta\colon \R \to \R$ be the Kronecker delta function, i.\ e.\ $\delta(0)=1$ and $\delta(x)=0$ otherwise (this is used in \ref{c1approx}). Set $f_0 = f$, $K_0 = K$, $\coD_0 = \coL_{-\epsilon}$ and $\delta_0 = \eta-\epsilon$. We inductively construct a sequence of Legendrian curves $f_j \colon M \to \coD$ and an increasing sequence of compact sets $K_j \subset \Int M$, $j \in \N$, satisfying $K_j \subset \Int K_{j+1}$ and $\cup_j K_j = \Int M$, such that:
		\begin{enumerate}[label=(\alph*$_j$)]
			\item 
			\label{Kapprox}
			$\|f_{j}-f_{j-1}\|_{0,K_{j-1}} < \varepsilon_j$, where $\varepsilon_j$ is defined below;
			\item 
			\label{Mapprox}
			$\|f_{j}-f_{j-1}\|_{0,M} < \sqrt{2\delta_{j-1}^2+2\delta_{j-1}/\kappa_{\coD_{j-1}}^\mmin}$;
			\item
			\label{bddist}
			$\dist_{f_{j}}(p_0) > j$;
			\item 
			\label{j-1proper}
			$f_{j}(M \backslash \Int K_{j-1}) \subset \coD \backslash \overline{\coD}_{j-1}$;
			\item 
			\label{jproper}
			$f_{j}(M \backslash \Int K_{j}) \subset \coD \backslash \overline{\coD}_{j}$;
			\item 
			\label{j-1interpolation}
			$f_{j}$ agrees with $f_{j-1}$ to order at least $m$ at every point $p \in \Lambda$;	
			\item
			\label{jimmersion}
			$f_j$ is an immersion, if $\dd f(p) \neq 0$ for every $p \in \Lambda$,
			\item 
			\label{jembedding}
			$f_j$ is an embedding, if the condition in \ref{jimmersion} holds and $f|_{\Lambda}$ is injective, and		
			\item
			\label{c1approx}
			$\varepsilon_j = \min\{\varsigma_j, \tau_j, \varepsilon/2^j\}$, where $\tau_j = \max\{\tau_j',\delta(\tau_j')\}$, $\varsigma_j = \max\{\varsigma_j', \delta(\varsigma_j')\}$, and
			\[
			\tau_j' = \frac{1}{2}\min_{p\in M} \| \partial f_{j-1} (p)\| \quad \text{and} \quad
			\varsigma_j' = \frac{1}{2j^2}\inf \left\{\|f_{j-1}(p)- f_{j-1}(q)\|: p,q \in M, \;
			\dist(p,q) > \frac{1}{j} \right\}.
			\]
		\end{enumerate}

		We claim that such a sequence converges to a holomorphic Legendrian curve $F\colon M \to \C^{2n+1}$ satisfying the required properties. Note that by \ref{Mapprox} and the definition of $\delta_j$, the maps $f_j$ indeed converge uniformly on $M$, hence the limit map $F$ is continuous and maps $M$ into $\coD$. Moreover, by \ref{j-1proper} and \ref{jproper} the map $F$ satisfies $F(M \backslash \Int K) \subset \coD \backslash \overline{\coL}$. By \ref{Kapprox}, the convergence is uniform on compacts in $\Int M$, thus the limit $F|_{\Int M}$ is holomorphic and Legendrian. By \ref{Kapprox} and \ref{bddist} the limit map satisfies $\dist_{F}(p_0) = \infty$ and $\dist_F(p) = \infty$ for any $p \in M$, hence the map $F$ is also complete. By \ref{Kapprox} and \ref{c1approx}, the limit map stays $\varepsilon$-close to $f$ in the $\c^1(K)$-topology. By \ref{j-1interpolation} the map $F$ agrees with $f$ to a order at least $m$ at every point $p \in \Lambda$. 
		
		Now suppose $\dd f(p) \neq 0$ for every $p \in \Lambda$. Then \ref{jimmersion} holds for every $j$, thus $\tau_j' > 0$. Any point $p \in \Int M$ is contained in some $K_{j_0}$, hence in all $K_j$ for $j \geq j_0$. 
		Then, for any such $j$ the following holds by \ref{Kapprox} and \ref{c1approx}:
		\begin{align*}
			\|\partial f_{j}(p)\| & \leq \|\partial f_{j+1}(p)\| + \|\partial f_{j+1}(p) - \partial f_j(p)\|
			< \\
			& <
			\|\partial f_{j+1}(p)\| + \varepsilon_{j+1} \leq 
			\|\partial f_{j+1}(p)\| + \frac{1}{2}\min_{p\in M} \|\partial f_j(p) \|,
		\end{align*}
		thus $ \| \partial f_{j+1}(p) \| > 2^{-1} \min_{p\in M} \|\partial f_j(p) \|$. By induction we have for any $k \in \N$ that $\| \partial f_{j+k}(p) \| > 2^{-k} \min_{p\in M} \|\partial f_j(p) \|$, hence	
		\begin{align*}
			\|\partial F(p)\| & \geq \|\partial f_{j_0}(p)\| - \|\partial F(p) - \partial f_{j_0}(p)\| \geq \\
			& \geq \|\partial f_{j_0}(p)\| - 
			\sum_{j=j_0}^{\infty} \|\partial f_{j+1}(p) - \partial f_j(p) \|
			\geq \\
			&\geq
			\|\partial f_{j_0}(p)\| - 
			\sum_{k=0}^{\infty} \frac{1}{2} \min_{q \in M} \|\partial f_{j_0+k}(q)\|
			\\
			& \geq \|\partial f_{j_0}(p)\| -  
			\sum_{k=0}^{\infty} \frac{1}{2^{k+1}}\min_{q\in M} \| \partial f_{j_0} (q)\| > 0,
		\end{align*}
		proving that $F|_{\Int M}$ is an immersion. 
		
		If, in addition, $f|_{\Lambda}$ is injective, then \ref{jembedding} holds for every $j$, hence $\varsigma_j' > 0$.	For a pair of distinct points $p,q \in \Int M$ find a $j_0 \in \N$ such that $d(p,q) > 1/j_0$. Then $d(p,q) > 1/j$ for every $j \geq j_0$ and we have by \ref{Kapprox} and \ref{c1approx}:
		\begin{align*}
			\|f_{j}(p)-f_{j}(q)\| &\leq \|f_{j+1}(p)-f_j(p)\| +
			\|f_{j+1}(q) - f_j(q)\| + \|f_{j+1}(p) - f_{j+1}(q)\| < \\
			&<
			2\varsigma_{j+1} + \|f_{j+1}(p) - f_{j+1}(q)\| \leq \\
			& \leq
			\frac{1}{(j+1)^2} \|f_{j}(p)-f_{j}(q)\| + \|f_{j+1}(p) - f_{j+1}(q)\| \\
			\|f_{j+1}(p) - f_{j+1}(q)\| &\geq \left(1-\frac{1}{(j+1)^2}\right) \|f_{j}(p)-f_{j}(q)\|.
		\end{align*}
		By induction we obtain the following for any $k \in \N$:
		\[
		\|f_{j_0+k}(p) - f_{j_0+k}(q)\| \geq 
		\|f_{j_0}(p) - f_{j_0}(q)\|
		\prod_{j=j_0}^{j_0+k}\left(1-\frac{1}{j^2}\right)
		\]
		The product on the left converges to a strictly positive value $C>0$ as $k \to \infty$, thus in the limit the following holds:
		\[
		\|F(p) - F(q)\| \geq C\|f_{j_0}(p) - f_{j_0}(q)\| > 0.
		\]
		Since $p,q$ were arbitrary, the map $F|_{\Int M}$ is injective. In addition, the properties \ref{j-1proper} and \ref{jproper} guarantee $F|_{\Int M}$ is proper. But, by the above, $F|_{\Int M}$ is also an immersion, hence a proper embedding.
		
		We now explain how to obtain the sequence $f_j$. Note that $f_0=f$ (defined above) satisfies (e$_0$) while (a$_0$), (b$_0$), (c$_0$), (d$_0$), and (f$_0$) are void. Conditions (g$_0$) and (h$_0$) are obtained by applying \cite[Theorem 2.1]{Svetina2024} to $f_0$ provided the suitable conditions hold.
		
		Suppose $f_1,\ldots,f_{j-1}$ as above have been constructed for some $j \in \N$. Apply Lemma \ref{cvxInductiveLemma} to the data $\coL = \coD_{j-1}$, $\coD = \coD$, $f=f_{j-1}$, $K = K_{j-1}$, $\eta = \delta_{j-1}$, $\delta=\delta_j$, to obtain a holomorphic Legendrian curve $F=f_j' \colon M \to \coD$ satisfying \ref{Kapprox}, \ref{Mapprox}, \ref{j-1proper}, \ref{j-1interpolation}, and \ref{c1approx}. Moreover, the lemma guarantees that $f_{j-1}'(bM) \subset \coD \backslash \overline{\coD}_{j}$. Then, use Lemma \ref{bdincrease} to approximate $f_j'$ with a curve $f_j$ satisfying \ref{bddist}, and \ref{j-1interpolation}. Again, use \cite[Theorem 2.1]{Svetina2024} to obtain \ref{jimmersion} or \ref{jembedding} if needed. If these last approximations are good enough, then $f_j(bM) \subset \coD \backslash \overline{\coD}_{j}$, hence there exists a compact set $K_j \subset \Int M$, such that $K_{j-1} \subset \Int K_j$ and $f(M\backslash \Int K_j) \subset \coD \backslash \overline{\coD}_j$, thus obtaining property \ref{jproper}. Obviously all the other properties still hold by choosing a well enough approximation of $f_{j-1}'$ by $f_j$, thus closing the induction.
	\end{proof}

	\begin{rmk}
		\label{globalMapprox}
		If $f(bM)$ is close enough to $b\coD$, meaning there exists an $\eta < 1/\kappa_{\coD}^\mmax$ such that $f(bM) \subset \coD \backslash \overline{\coD}_\eta$, then by \ref{deltainequality}, the limit curve $F$ in this case satisfies
		\[
		\|F-f\|_{0,M} < \sqrt{2\eta^2+ \frac{2\eta}{\kappa_\coL^\mmin}} = \o(\sqrt{\eta}), \quad
		\text{as} \quad \eta \to 0.
		\]
	\end{rmk}


	Our second theorem generalises Theorem \ref{properCvx} to arbitrary convex domains at the cost of losing control over the uniform norm of the difference between the approximating mapping and the approximant. As mentioned in the introduction, Theorem \ref{weakCvx} is an immediate consequence of the following result.
	
	\begin{thm}
		\label{convExhaustionThm}
		Suppose $\coD \subset \C^{2n+1}$ is a convex domain in $\C^{2n+1}$ and let $f \colon M \to \coD$ be a holomorphic Legendrian curve, defined on a compact bordered Riemann surface $M$. Given a compact set $K \subset \Int M$, a finite set $\Lambda' \subset \Int K$, a number $m \in \N\cup \{0\}$, a closed discrete set $\Lambda'' \subset \coD\backslash f(M)$, and a number $\varepsilon > 0$, there exists a proper and complete holomorphic Legendrian curve $F\colon \Int M \to \coD$ such that
		\begin{enumerate}[label=\Roman*)]
			\item $\|F-f\|_{0,K} < \varepsilon$,
			
			\item $F$ agrees with $f$ to order at least $m$ at every point $p \in \Lambda'$,
			
			\item $\Lambda'' \subset F(\Int M)$.
		\end{enumerate}
		If $m=0$ or $m\geq 1$ and $\dd f(p) \neq 0$ for every point $p \in \Lambda'$, then $F$ may be taken an immersion. If, furthermore, $f|_{\Lambda'}$ is injective, then $F$ may be made an embedding.
	\end{thm}

	\begin{proof}
		By \cite[Lemma 2.3.2]{Hoermander2007}, any bounded convex domain admits an exhaustion $\coD_1 \Subset \coD_2 \Subset \cdots$ by smoothly bounded strongly convex compact sets $\coD_j$, $j \in \N$. By the same lemma, this also holds for unbounded sets. Namely, for such a set $\coD$ define $\coD_j' := B(0,j) \cap \{x \in \coD: \dist(x,b\coD)>1/j\}$ and note that $\coD_j'$ is a bounded convex domain. Then approximate $\coD_j'$ with smoothly bounded strongly convex sets $\coD_j$ such that $\coD_j' \subset \coD_j \subset K_{j+1}'$. Since $\Lambda''$ is closed and discrete, we may assume $\Lambda''$ does not intersect $b\coD_j$ for any $j \in \N$.
		
		Set $\Lambda_j'' = \Lambda'' \cap \coD_j\backslash \overline{\coD}_{j-1}$ and $\coD_0 = \emptyset$, $f_0 := f$. We will construct a sequence of holomorphic Legendrian curves $f_j \colon M \to \coD$, an increasing sequence $\{K_j\}_{j\in \N}$ of compact sets $K_j \subset \Int M$ satisfying $K_j \subset \Int K_{j+1}$ for every $j \in \N$ and $\cup_j K_j = \Int M$, and an increasing sequence of finite sets $\Lambda_j' \subset \Int M$, $\Lambda' \subset \Lambda_j$, satisfying $\Lambda_j''\subset f_j(\Lambda_j')$ such that:
		\begin{enumerate}[label=(\alph*$_j$)]
			\item 
			\label{Kapprox'}
			$\|f_{j}-f_{j-1}\|_{1,K_{j-1}} < \varepsilon_j$, where $\varepsilon_j$ is defined below;
			\item
			\label{bddist'}
			$\dist_{f_{j}}(p_0, bK_j) > j$;
			\item 
			\label{j-1proper'}
			$f_{j}(M \backslash \Int K_{j-1}) \subset \coD \backslash \overline{\coD}_{j-1}$;
			\item 
			\label{jproper'}
			$f_{j}(M \backslash \Int K_{j}) \subset \coD \backslash \overline{\coD}_{j}$;
			\item 
			\label{j-1interpolation'}
			$f_{j}$ agrees with $f_{j-1}$ to order at least $m$ at every point $p \in \Lambda_{j-1}'$;
			\item 
			\label{jinterpolation'}
			$\Lambda''_j \subset f_j(\Int M)$,
			
			\item
			\label{jimmersion'}
			$f_j$ is an immersion, if $\dd f(p) \neq 0$ for every $p \in \Lambda'$,
			\item 
			\label{jembedding'}
			$f_j$ is an embedding, if the condition in \ref{jimmersion'} holds and $f|_{\Lambda'}$ is injective, and
			
			\item 
			\label{c1approx'}
			$\varepsilon_j = \min\{\varsigma_j, \tau_j, \varepsilon/2^j\}$, where $\tau_j = \max\{\tau_j',\delta(\tau_j')\}$, $\varsigma_j = \max\{\varsigma_j', \delta(\varsigma_j')\}$, and
			\[
			\tau_j' = \frac{1}{2}\min_{p\in M} \| \partial f_{j-1} (p)\| \quad \text{and} \quad
			\varsigma_j' = \frac{1}{2j^2}\inf \left\{\|f_{j-1}(p)- f_{j-1}(q)\|: p,q \in M, \;
			\dist(p,q) > \frac{1}{j} \right\}.
			\]
		\end{enumerate}
		
		Suppose such a sequence exists. By \ref{Kapprox'}, the sequence converges uniformly on compact sets in $\Int M$ to a holomorphic Legendrian curve $F \colon \Int M \to \coD$ satisfying $\|F-f\|_K < \varepsilon$. By \ref{j-1proper'} and \ref{jproper'}, it follows that the curve $F$ is proper. By \ref{j-1interpolation'} and \ref{jinterpolation'} the curve $F$ satisfies the required interpolation condition and by \ref{bddist'} the obtained limit curve $F$ is complete. Meanwhile, conditions \ref{jimmersion'}, \ref{jembedding'} and \ref{c1approx'} imply the limiting curve is an immersion or an embedding provided that suitable conditions in the theorem hold. This is checked in the same way as in the proof of Theorem \ref{properCvx}.
		
		Let us now show how to construct such a sequence. Suppose $f_1,\ldots, f_{j-1}$, and $K_1, \ldots K_{j-1}$, have already been constructed. Use Lemma \ref{inductionStart} to approximate $f_{j-1}$ uniformly on $K_{j-1}$ by a holomorphic Legendrian curve $f_{j-1}'\colon M \to \coD$ such that $f_{j-1}'(M \backslash \Int K_{j-1}) \subset \coD \backslash \overline{\coD}_{j-1}$ and $\Lambda_j'' \subset f_{j-1}'(\Int M)$. For every $p \in \Lambda_j''$ choose a unique point $u_p \in \Int M$ satisfying $f_{j-1}'(u_p) = p$ and let $\Lambda_j' = \Lambda_{j-1}' \cup  \{u_p\}_{p \in \Lambda_j''}$. For every point $p \in \Lambda_j'$ choose a small embedded disk $\Omega_p \subset M$ with center $p$ and let $K_{j-1}' = K_{j-1} \cup \{\overline{\Omega}_p\}_{p \in \Lambda_j'}$. Then, choose a small $t > 0$ and let $\coD'_j = (\coD_j)_{-t}$ in the sense of \ref{mt}, that is, we slighlty enlarge $\coD_j$ by pushing its boundary in the direction of the \emph{outer} unit normal vector field for time $t$. Use Theorem \ref{properCvx} to approximate $f_{j-1}'$ on $K_{j-1}'$ by a proper and complete holomorphic Legendrian curve $f_{j}\colon M \to \coD_j'$ which agrees with $f_{j-1}'$ to order $m$ at every point in $\Lambda_j' \subset \Int K_{j-1}'$, thus $f_j$ satisfies properties \ref{Kapprox'}, \ref{bddist'}, \ref{j-1interpolation'}, and \ref{jinterpolation'}. Moreover, we ensure that $f_{j}(M \backslash \Int K_{j-1}) \subset \coD \backslash \overline{\coD}_{j-1}$, thus satisfying property \ref{j-1proper'}. Since $f_{j}(bM) \subset \coD \backslash \overline{\coD}'_j \subset \coD \backslash \overline{\coD}_j$, there exists a compact set $K_j \subset \Int M$ such that $K_{j-1} \cup \Lambda_j' \subset \Int M$ and $f_j(M \backslash \Int K_j) \subset \coD \backslash \overline{\coD}_j$, thus property \ref{jproper'} holds. Using \cite[Theorem 3.5]{Svetina2024} we may ensure that $f_j$ is an immersion or an embedding, thus satisfying \ref{jimmersion'} and/or \ref{jembedding'}, while still satisfying property \ref{bddist'}, closing the induction.	
	\end{proof}
	
	\section{Approximation of curves hitting the boundary}
	\label{properWcomplete}
	In this section we consider the problem of approximation of holomorphic Legendrian curves in closed convex domains in Euclidean spaces whose image does not lie in the interior of the domain but may intersect the boundary somewhere, compare with part b) in \cite[Theorem 1.2]{Alarcon2015} for the similar problem in the case of minimal surfaces. In the cited paper the authors use a simple homothety to perturb the given minimal surface and reduce case b) to the case where entire image of the compact Riemann surface is contained in the interior of the convex domain. This is possible by the fact that the property of being a minimal surface is invariant under rigid motions and homotheties of the Euclidean space. The same is not true in the setting of contact geometry, i.e., translations, rotations and homotheties do not preserve Legendrian curves. Instead, the role of translations is played by suitable \emph{contactomorphisms}, that is biholomorphic maps preserving the contact structure, while homotheties are replaced by flows of certain vector fields.
	
	Namely, a vector field $V$ on a complex contact manifold $(X,\xi)$ is called \emph{contact} if its flow preserves the contact structure meaning that $(\phi_t)_* \xi = \xi$ holds (locally) for all $t$ where the flow $\phi_t$ of $V$ is defined. If $\xi = \ker \alpha$ for a holomorphic 1-form $\alpha$, then $V$ is contact if and only if $\phi_t^*\alpha = \lambda \alpha$ for a holomorphic function $\lambda \in \o(X)$. It is well known that holomorphic contact vector fields on $X$ are in bijective correspondence with holomorphic functions on $X$ given by
	\begin{equation}
		\label{cntvf}
	\o(X) \ni f \mapsto V_f, \quad \text{such that }
	\alpha(V_f) = f, \text{ and }
	\iota_{V_f}\dd \alpha = -\dd f +R(f)\alpha,
	\end{equation}
	where $R$ is the \emph{Reeb vector field} of the contact form $\alpha$, see \cite[Theorem A.3]{Alarcon2017} and also \cite[Section 2.3]{Geiges2008}. $R$ is itself determined by the system of equations
	\[
	\alpha(R) \equiv 1, \quad \text{and} \quad \iota_R \dd \alpha = 0.
	\]

	\begin{defi}
		A convex domain $\coD \subset \C^{2n+1}$ is \emph{Legendrian convex}, if for any generalised Legendrian curve $f\colon M \to \overline{\coD}$, defined on a compact bordered Riemann surface $M$, there exists a neighbourhood $U \subset \C^{2n+1}$ of $f(M)$, a holomorphic contact vector field $V$, defined on $U$, whose flow is defined on $f(M) \times (-\varepsilon,\varepsilon)$ for some $\varepsilon > 0$, satisfying
		\[
		\phi_t(f(M)) \subset \coD,
		\]
		for every $t \in (0,\varepsilon)$. We call such $V$ an \emph{inward pointing contact vector field} for $M$ on $\overline{\coD}$.
	\end{defi}
	In particular, a relatively compact convex domain $\coD\subset \C^{2n+1}$ is Legendrian convex if the boundary of $\coD$ is of class $\c^1$ and there exists a holomorphic contact vector field $V$, defined on a neighbourhood of $\overline{\coD}$, that is transverse to $b\coD$ along the entire boundary $b\coD$. Compare this to the case of convex surfaces in real contact manifolds of dimension 3. Namely, a (hyper)surface $S \subset M$ in a real contact manifold $(M,\xi)$ is called \emph{convex} if there exists a contact vector field $V$ on a neighbourhood of $S$ that is everywhere transverse to $S$, see \cite[Section 4.6.2]{Geiges2008}.	
	
	In Legendrian convex domains we are able to push Legendrian curves that intersect the boundary into the interior using flows of inward pointing contact vector fields. If $f\colon M \to \overline{\coD}$ is such a Legendrian curve, then $\phi_t(f(M)) \subset \coD$ for some $t > 0$ and $(\phi_t \circ f)^* \alpha = f^*(\lambda \alpha) = 0$ for some holomorphic function $\lambda$ on a neighbourhood of $\overline{\coD}$, thus $\phi_t \circ f$ is a Legendrian curve whose image lies in the interior $\coD$ of $\overline{\coD}$.
	
	\vspace{10pt}
		The following lemma provides contact analogues of translations and homotheties.
	
	\begin{lm}
		The following transformations of $\C^{2n+1}$ preserve Legendrian curves:
		\begin{itemize}
			\item[(T)] \label{T} 
			$T_p\colon (\bx,\by,z) \mapsto (\bx+\bx_0, \by+\by_0, z+z_0-\bx_0\by)$ for any $p:=(\bx_0, \by_0, z_0) \in \C^{2n+1}$,
			
			\item[(H)] \label{H}
			$H_{\lambda} \colon (\bx,\by,z) \mapsto (\lambda \bx, \lambda \by, \lambda^2 z)$ for any $\lambda \in \C^*$.
		\end{itemize}
	\end{lm}
	
	\begin{proof}
		It follows by a short computation that $T_p^*\alpha = \alpha$ and $H_p^*\alpha = \lambda^{2}\alpha$ for every $p \in \C^{2n+1}$. Hence, if $f\colon M \to \C^{2n+1}$ is a Legendrian curve, then $(T \circ f)^*\alpha = f^*(T^*\alpha) = f^*\alpha= 0$ and similarly $(H \circ f)^*\alpha = \lambda^2 f^*\alpha = 0$.
	\end{proof}
	
	Since both $H$ and $T$ are linear, they also preserve convexity. However, not every convex domain is closed under the action of maps $\{H_\lambda\}_{\lambda \in [0,1]}$. For example, take 
	\[
	\coD = \{(x,y,z) \in \C^3: \Re z > \Re x\}
	\]
	and let $p=(x,y,z) \in b\coD$ be a boundary point satisfying $\Re z > 0$. Then, for any $\lambda < 1$, we have
	\[
	\Re \pr_3(H_\lambda(p)) = \lambda^2 \Re z = \lambda^2 \Re x < \lambda \Re x = \pr_1(H_\lambda(p)),
	\]
	thus $p$ leaves $\overline{\coD}$ upon the action of $H_\lambda$ for any $\lambda < 1$.
	
	Note that by replacing $\lambda$ with $e^{t}$ and taking the $t$-derivative of $H_\lambda$ we obtain a holomorphic vector field
	\begin{equation}
		\label{inwCont}
	V(\bx,\by,z) = 
	\dv{t} \left(\bx_0 e^t, \by_0 e^t, z_0e^{2t} \right) =
	\left(\bx, \by, 2z\right)
	\end{equation}	
	Hence, if $\coD=\{\rho = r^2\}$, $\rho(p) = \|p\|^2$, is a ball of radius $r$ in $\C^{2n+1}$, centered at the origin, then $\coD$ is Legendrian convex, since 
	\[
	\frac{1}{2}\dd \rho_p(V_p) = - \|\bx\|^2-\|\by\|^2-2|z|^2 < 0.
	\]
	Consequently, every ball $B(0,r)$ is Legendrian convex. 
	
	On the other hand, let $P = r_1\D \times \cdots \times r_{2n+1}\D \subset \C^{2n+1}$ be a polydisk of polyradius $R = (r_1,\ldots,r_{2n+1})$, $r_i > 0$ for every $i=1,\ldots,2n+1$, and let $V= -(\bx,\by,2z)$ be the negative of the holomorphic contact vector field \ref{inwCont}. Let $\pr_i\colon \C^{2n+1} \to \C$ denote the projection onto the $i$-th coordinate line and write
	\[
	P = \bigcap_{i=1}^{2n+1} \{|\pr_i|<r_i\}.
	\]
	Hence, if $p=(\bx,\by,z) \in \overline{P}$, then $\phi_t(p) = (e^{-t}\bx, e^{-t}\by, e^{-2t}z) \in P$, thus every such polydisk is Legendrian convex.

	\begin{thm}
		Suppose $\coD$ is a Legendrian convex domain in $\C^{2n+1}$ and $f\colon M \to \overline{\coD}$ is a holomorphic Legendrian curve, defined on a compact bordered Riemann surface $M$.
		\begin{enumerate}[label=\Roman*)]
			\item 
			\label{somewhereProper}
			The curve $f$ may be approximated uniformly on compacts in $\Int M$ by continuous maps $F \colon M \to \overline{\coD}$ such that $F|_{\Int M}\colon \Int M \to \coD$ is a proper and complete holomorphic Legendrian embedding. Moreover, for a given discrete set $\Lambda \subset \Int M$, one may find $F$ such that $\Lambda \subset F(\Int M)$.
			
			\item
			\label{alreadyproper}
			If $f(bM) \subset b\coD$, then $f$ may be approximated uniformly on $M$ by continuous maps $F \colon M \to \overline{\coD}$ such that $F|_{\Int M}\colon \Int M \to \coD$ is a proper and complete holomorphic Legendrian embedding.
		\end{enumerate}
	\end{thm}
	
	\begin{rmk}
		Note that the only novelty in case \ref{somewhereProper} is when the starting curve $f\colon M \to \overline{\coD}$ hits the boundary somewhere, since if $f(bM) \subset \coD$ this is just Theorem \ref{convExhaustionThm}. In case \ref{alreadyproper}, a hitting result similar to the one in \ref{somewhereProper} is not possible as one wishes to perturb $f$ as little as possible in order to obtain uniform approximation on the entire surface $M$. On the other hand, uniform approximation on $M$ is not possible in case \ref{somewhereProper}, as one does not know a priori how far from the boundary of $b\coD$ may $f$ take some points in $bM$. Note also that the theorem does not include interpolation of the starting curve as this is ruined by deforming its image with the above versions of homotheties and translations. Of course, in \ref{somewhereProper} one does not need the inward pointing contact vector field to be defined along the entire boundary $b\coD$ but only at the intersection $f(bM) \cap b\coD$.
	\end{rmk}
	
	\begin{proof}
		So suppose $f(bM) \cap b\coD$ is nonempty. Take an inward pointing holomorphic contact vector field $V$ for $M$ $\overline{\coD}$ and choose a suitably small $t>0$ such that $\phi_t(f(M)) \subset \coD$ and $\phi_t \circ f$ is sufficiently uniformly close to $f$ on $M$. Then use Theorem \ref{convExhaustionThm} with the data $K$, $\Lambda'' = \Lambda$, and arbitrary $\varepsilon > 0$ to approximate $\phi_t \circ f$ uniformly on $K$ by a continuous map $F\colon M \to \overline{\coD}$ such that $F|_{\Int M}$ is a proper and complete Legendrian embedding and $\Lambda \subset F(\Int M)$, proving \ref{somewhereProper}.
		
		To prove \ref{alreadyproper} note that $f(bM) \subset b\coD$ and a choice of $t>0$ sufficiently close to 0 imply $\phi_t(f(bM)) \subset \coD \backslash \overline{\coD}_\eta$ for some $\eta > 0$ which can be made arbitrarily small by choosing $t$ close enough to 0. By Remark \ref{globalMapprox} we see that $F$ constructed as above satisfies 
		\[
		\|F-\phi_t \circ f\|_{0,M} < \sqrt{2\eta^2+2\eta/\kappa_{\coD}^\mmin} = \o(\sqrt{\eta}),
		\]
		completing the proof.
	\end{proof}

	\begin{proof}[Proof of Corollary \ref{properExst}]
		By choosing an interior point $p_0 \in \coD$ and replacing $\coD$ with $T_{p_0}^{-1}\coD$ we may assume $\coD$ contains the origin in $\C^{2n+1}$. Choose an arbitrary holomorphic Legendrian curve $\tilde{f}\colon M \to \C^{2n+1}$ which exists by \cite[Theorem 1.1]{Alarcon2017}. Choosing $\lambda > 0$ small enough, the curve $H_\lambda \circ \tilde{f}$ can be made to lie in $\coD$ and miss every point in $\Lambda''$, since the latter set is discrete. Now use Theorem \ref{convExhaustionThm} with the data $K, \Lambda' = \emptyset$, $\Lambda'' = \Lambda$ and arbitrary $\varepsilon > 0$ to obtain a holomorphic Legendrian curve $F \colon M \to \coD$ satisfying the required properties.
	\end{proof}

	\section{Almost proper curves}
	\label{almostProperProof}
	
	In this section we prove Theorem \ref{denseImmersion} as a special consequence of the following result.
	
	\begin{thm}
		\label{almostProperApprox}
		Suppose $f\colon M \to \coD$ is a holomorphic Legendrian curve, defined on a compact bordered Riemann surface $M$, taking values in a convex domain $\coD \subset \C^{2n+1}$. Given a compact set $K \subset \Int M$, a finite set $\Lambda' \subset K$, a number $m \in \N \cup \{0\}$, a countable set $\Lambda'' \in \coD\backslash f(\Lambda')$, and a number $\varepsilon > 0$, there exists an almost proper and complete holomorphic Legendrian curve $F \colon \Int M \to \coD$, such that
		\begin{enumerate}[label=\Roman*)]
			\item $\|F-f\|_{0,K} < \varepsilon$,
			
			\item $F$ agrees with $f$ to order at least $m$ at every point $p \in \Lambda'$, and
			
			\item $\Lambda'' \subset F(\Int M)$.
		\end{enumerate}
	\end{thm}
	
	If $\dd f(p) \neq 0$ for every $p \in \Lambda'$, then $F$ may be made an immersion. If, additionaly, the map $f|_{\Lambda'}$ is injective, then $F$ may be made injective.
	
	\begin{proof}
		Fix an enumeration $\lambda_1, \lambda_2, \ldots$ for the set $\Lambda''$, let $\Lambda_j'' := \{\lambda_1,\ldots,\lambda_j\}$ for every $j \in \N$, $\Lambda_0'= \Lambda'$, fix a normal exhaustion $K_0 := K \Subset K_1 \Subset K_2 \Subset \cdots$ of $\Int M$ by Runge compact sets $K_j$ and an exhaustion $\coD_0 \Subset \coD_1 \Subset \cdots$ of $\coD$ by smoothly bounded strongly convex relatively compact sets $\coD_j \Subset \coD$. Set $f_0 := f \colon M \to \coD$. We will construct a sequence $f_j \colon M \to \coD$ of holomorphic Legendrian curves satisfying the following properties:
		\begin{enumerate}[label=\alph*$_j$)]
			\item $\|f_j - f_{j-1}\|_{1,K_{j-1}} < \varepsilon_j$ where $\varepsilon_j > 0$ will be defined later on,
			
			\item $d_{f_j}(p_0) > j$,
			
			\item $f_j(bK_j) \subset \coD \backslash \overline{\coD}_j$,
			
			\item $f_j$ agrees with $f_{j-1}$ to order at least $m$ at every point $p \in \Lambda_j'$, where the set $\Lambda_j'$ will be defined later on,
			
			\item $\Lambda_j''\subset f_j(\Int M)$,
			
			\item $f_j$ is an immersion, if $\dd f(p) \neq 0$ for every $p \in \Lambda'$,
			
			\item $f_j$ is an embedding, if also $f|_{\Lambda'}$ is injective, and
			
			\item define $\tau_j$, $\varsigma_j$ as in the proof of Theorem \ref{convExhaustionThm}, let
			\[
			\rho_j := 2^{-j} \min_{l=1,\ldots,j-1} \dist(f_{l}(bK_{l}), b\coD_{l}) 
			\]
			and let $\varepsilon_j = \min \{\varsigma_j, \tau_j, \rho_j, 2^{-j}\varepsilon\}$.
		\end{enumerate}
		
		Suppose such a sequence exists. Provided $\sum_j\varepsilon_j$ converges, the properties a$_j$) guarantee that the sequence of curves $f_j$ converges uniformly on compacts in $\Int M$ to a holomorphic Legendrian curve $F\colon M \to \coD$. The properties b$_j$) guarantee that $F$ is complete. Properties d$_j$) guarantee that $F$ agrees with $f$ to order $m$ at every point $p \in \Lambda'$ and, together with properties e$_j$), that $\Lambda'' \subset F(\Int M)$. Now suppose $L \subset \coD$ is compact and let $j_0 \in \N$ be such that $L \subset \coD_{j_0}$, hence $L \subset \coD_j$ for every $j \geq j_0$. Choose a connected component $C \subset F^{-1}(L)$ and suppose that $C \cap \Int K_{j}$ is nonempty for some $j \geq j_0$. By c$_{j}$) we obtain $f_{j}(bK_{j}) \subset \coD \backslash \overline{\coD}_j$ and by noting that the properties h$_j$) imply $\varepsilon_l \leq 2^{-l} \dist(f_{j}(bK_{j}), b\coD_{j})$ for every $l > j$ the following holds by a$_j$) for every $p \in bK_j$:
		\[
		\|F(p)-f_j(p)\| \leq 
		\sum_{l > j}\|f_{l}(p)-f_{l-1}(p)\| <
		\sum_{l > j} 2^{-l} \dist(f_{j}(bK_{j}), b\coD_{j}) =
		2^{-j} \dist(f_{j}(bK_{j}), b\coD_{j}).
		\]
		Since $f_j(p) \in \coD\backslash \overline{\coD}_j$, this implies $F(p) \in \coD\backslash \overline{\coD}_j$. Thus, $C$ contains no boundary point of $K_j$, hence $C \subset K_j$, i.~e. $C$ is compact. That $F$ is an immersion or an injective immersion follows from the properties g$_j$), -- i$_j$), and is checked in the same way as in the proof of Theorem \ref{convExhaustionThm}.
		
		Now suppose $f_1,\ldots,f_{j-1}$ have been constructed. Assume $M$ is a smoothly bounded compact domain in an open Riemann surface $\widetilde{M}$. Choose a smooth Jordan arc $\gamma_j$ in $\widetilde{M}$ with one endpoint in a boundary point $p \in bM$ such that $\gamma_j$ intersects $M$ only in the point $p$ and such that this intersection is transverse. We label the other endpoint of $\gamma_j$ by $q \in \widetilde{M} \backslash M$. Then, extend $f_{j-1}$ to a generalised Legendrian curve $f_{j-1}'\colon M \cup \gamma_n \to \coD$ such that $f_{j-1}'(q') = \lambda_n$ for some point $q'$ in the interior of the arc $\gamma_n$, and the curve $f_{j-1}'$ is holomorphic near $q'$ with $\dd f_{j-1}'(q') \neq 0$. By \cite[Theorem 2.1]{Svetina2024} there exists a holomorphic Legendrian curve $f_{j-1}''\colon U \to \C^{2n+1}$, defined on a neighbourhood $U$ of $M$ in $\widetilde{M}$, approximating $f_{j-1}'$ on $K_{j-1}$, such that $f_{j-1}''$ agrees with $f_{j-1}'$ to order $m$ at every point in $\Lambda_j' := \Lambda_{j-1}'\cup\{q'\}$. In particular, $f_{j-1}''$ agrees with $f_{j-1}$ to order at least $m$ at every point in $\Lambda'$ and $\Lambda_j'' \subset f_{j-1}''(\Int M)$. By shrinking $U$ if necessary, we may assume $f_{j-1}''(M) \subset \coD$. 
		
		By Theorem \ref{fwbdexpose} there exists a conformal diffeomorphism $\phi\colon M \to \overline{U}$, holomorphic on the interior of $M$, which agrees with the identity map to order $m$ at every point $p \in \Lambda_j'$ and lies arbitrarily close to the identity map on the subset $K_{j-1} \subset M$. Thus, the map $f_{j-1}'' \circ \phi$ approximates well the map $f_{j-1}$ on $K_{j-1}$.	Now, use Lemma \ref{cvxInductiveLemma} to approximate $f_{j-1}''\circ \phi$ by a holomorphic Legendrian curve $f_j\colon M \to \coD$ such that the properties a$_j$) -- h$_j$) are satisfied, thus closing the induction. 
	\end{proof}
	
	\begin{proof}[Proof of Theorem \ref{denseImmersion}]
		View $M$ as the interior of the compact bordered Riemann surface $\overline{M}$ and choose any holomorphic Legendrian curve $f\colon \overline{M} \to \C^{2n+1}$, given by \cite[Theorem 1.1]{Alarcon2017}. Then, choose a compact set $K \subset M$ and a countable dense subset $\Lambda \subset \coD$. Finally, use Theorem \ref{almostProperApprox} with the data $K$, $\Lambda' = \emptyset$, $\Lambda'' = \Lambda$ and any $\varepsilon > 0$, noting that the last assumptions in the theorem are void.
	\end{proof}

	\bigskip
	\noindent
	\textsc{Acknowledgements:} The author is supported by grant MR-54828 from ARIS, Republic of Slovenia, associated to the research program P1-0291 \emph{Analysis and Geometry}. The author would like to thank F. Forstnerič for his advice, support and guidance.

	\bibliographystyle{abbrv}	
	\bibliography{document}

\begin{thebibliography}{10}

\bibitem{Alarcon2015}
A.~Alarc\'{o}n, B.~Drinovec~Drnov\v{s}ek, F.~Forstneri\v{c}, and F.~J.
  L\'{o}pez.
\newblock Every bordered {R}iemann surface is a complete conformal minimal
  surface bounded by {J}ordan curves.
\newblock {\em Proc. Lond. Math. Soc. (3)}, 111(4):851--886, 2015.

\bibitem{Alarcon2023}
A.~Alarcon and F.~Forstneric.
\newblock Oka-1 manifolds.
\newblock {\em arXiv e-prints}, page arXiv:2303.15855, Mar. 2023.

\bibitem{Alarcon2017}
A.~Alarc\'{o}n, F.~Forstneri\v{c}, and F.~J. L\'{o}pez.
\newblock Holomorphic {L}egendrian curves.
\newblock {\em Compos. Math.}, 153(9):1945--1986, 2017.

\bibitem{mincplx}
A.~Alarc\'{o}n, F.~Forstneri\v{c}, and F.~J. L\'{o}pez.
\newblock {\em Minimal surfaces from a complex analytic viewpoint}.
\newblock Springer Monographs in Mathematics. Springer, Cham, 2021.

\bibitem{Bellettini2013}
G.~Bellettini.
\newblock {\em Lecture notes on mean curvature flow, barriers and singular
  perturbations}, volume~12 of {\em Appunti. Scuola Normale Superiore di Pisa
  (Nuova Serie)}.
\newblock Edizioni della Normale, Pisa, 2013.

\bibitem{Forstneric2017a}
F.~Forstneri\v{c}.
\newblock {\em Stein manifolds and holomorphic mappings}, volume~56 of {\em
  Ergebnisse der Mathematik und ihrer Grenzgebiete. 3. Folge. A Series of
  Modern Surveys in Mathematics}.
\newblock Springer, Cham, second edition, 2017.
\newblock The homotopy principle in complex analysis.

\bibitem{Forstneric2009}
F.~Forstneri\v{c} and E.~F. Wold.
\newblock Bordered {R}iemann surfaces in {$\mathbb{C}^2$}.
\newblock {\em J. Math. Pures Appl. (9)}, 91(1):100--114, 2009.

\bibitem{Geiges2008}
H.~Geiges.
\newblock {\em An introduction to contact topology}, volume 109 of {\em
  Cambridge Studies in Advanced Mathematics}.
\newblock Cambridge University Press, Cambridge, 2008.

\bibitem{Goluzin1969}
G.~M. Goluzin.
\newblock {\em Geometric theory of functions of a complex variable}, volume
  Vol. 26 of {\em Translations of Mathematical Monographs}.
\newblock American Mathematical Society, Providence, RI, 1969.

\bibitem{Gunning1967}
R.~C. Gunning and R.~Narasimhan.
\newblock Immersion of open {R}iemann surfaces.
\newblock {\em Math. Ann.}, 174:103--108, 1967.

\bibitem{Hoermander2007}
L.~H\"{o}rmander.
\newblock {\em Notions of convexity}.
\newblock Modern Birkh\"{a}user Classics. Birkh\"{a}user Boston, Inc., Boston,
  MA, 2007.
\newblock Reprint of the 1994 edition [of MR1301332].

\bibitem{Stout1965}
E.~L. Stout.
\newblock Bounded holomorphic functions on finite {R}eimann surfaces.
\newblock {\em Trans. Amer. Math. Soc.}, 120:255--285, 1965.

\bibitem{Svetina2024}
A.~Svetina.
\newblock Approximation of holomorphic {L}egendrian curves with
  jet-interpolation.
\newblock {\em J. Math. Anal. Appl.}, 531(2):Paper No. 127839, 25, 2024.

\end{thebibliography}

\end{document}